\documentclass[hidelinks,3p,10pt]{elsarticle}
\usepackage[utf8]{inputenc}
\usepackage{amssymb}
\usepackage{latexsym}
\usepackage{amsmath}
\usepackage{array}
\usepackage{graphicx}
\graphicspath{Images/}
\usepackage{caption}
\usepackage{subcaption}
\usepackage{algorithmic}
\usepackage{float}
\usepackage{hyperref}
\usepackage{cleveref}
\usepackage{algorithm}
\usepackage{amsthm}
\usepackage{verbatim}
\usepackage{color}

\newtheorem{theorem}{Theorem}

\newtheorem{lemma}{Lemma}
\newtheorem{remark}{Remark}

\newtheorem{corollary}{Corollary}

\journal{Applied Mathematics and Computation}

\begin{document}

\begin{frontmatter}

\title{Entropy conserving/stable schemes for a vector-kinetic model of hyperbolic systems}
\author {Megala Anandan\corref{cor1}}
\ead{megalaa@iisc.ac.in} 
\author {S.V. Raghurama Rao} 
\ead{raghu@iisc.ac.in}
\cortext[cor1]{Corresponding author}
\address{Indian Institute of Science, C.V. Raman Road, 560012, Bangalore, India} 
\date{October 26, 2023}

\begin{abstract}
The moment of entropy equation for vector-BGK model results in the entropy equation for macroscopic model. However, this is usually not the case in numerical methods because the current literature consists mostly of entropy conserving/stable schemes for macroscopic model. In this paper, we attempt to fill this gap by developing an entropy conserving scheme for vector-kinetic model, and we show that the moment of this results in an entropy conserving scheme for macroscopic model. With the numerical viscosity of entropy conserving scheme as reference, the entropy stable scheme for vector-kinetic model is developed in the spirit of Tadmor \cite{10.2307/2008251}. We show that the moment of this scheme results in an entropy stable scheme for macroscopic model. The schemes are validated on several benchmark test problems for scalar and shallow water equations, and conservation/stability of both kinetic and macroscopic entropies are presented.  
\end{abstract}


\end{frontmatter}

\section{Introduction}
The connection between entropy functions and symmetrisability of hyperbolic systems was explained in \cite{HARTEN1983151,Hughes1986ANF}, and this led to entropy-based non-linear stability analysis of numerical schemes. In the seminal work in \cite{10.2307/2008251,tadmor_2003}, a general condition to conserve/dissipate entropy of a semi-discrete scheme for hyperbolic system was introduced.  Following this, many developments on fluxes satisfying entropy conservation/dissipation condition for various hyperbolic systems were made. These include developments specific for shallow water equations \cite{GASSNER2016291,WINTERMEYER2017200,https://doi.org/10.1002/fld.4766}, Euler's equations \cite{Barth,ISMAIL20095410,puppo_semplice_2011,chandrashekar_2013,Ray2013335,ray_chandrashekar_fjordholm_mishra_2016,GASSNER201639,CREAN2018410,Chizari2021,Yan2023}, Navier-Stokes equations \cite{YAMALEEV2019108897,MANZANERO2020109241,Ranocha2020} and magneto hydro-dynamics equations \cite{doi:10.1137/15M1013626}. Recently, several interesting studies such as, entropy stability for conservation laws with non-convex flux functions \cite{LeFloch2020}, and characterisation of stability \cite{Gassner2022} and robustness (for under-resolved flows) \cite{JChan2022} of high order entropy stable schemes were carried out. \\
On the other hand, kinetic entropy formulations were introduced for hyperbolic equations like multi-dimensional scalar conservation laws, isentropic Euler and full Euler equations \cite{Perthame1991,10.2307/2152725,Lions1994,Deshpande2583}. Discrete kinetic models with entropy considerations were also proposed for hyperbolic systems \cite{10.2307/2587356,NATALINI1998292,Bouchut1999,Bouchut2003,Berthelin2002RelaxationTI,SMAI-JCM_2018__4__1_0}. Specifically, in \cite{Bouchut1999} it was shown that the entropy inequalities for a hyperbolic system can be derived as minimisation of entropies of vector-kinetic equation with BGK model. This approach of obtaining entropy inequalities from kinetic-BGK models is a promising strategy to characterise weak solutions of hyperbolic systems \cite{PS2011}.  Hence, in this paper, we attempt to develop entropy stable schemes  (in the sense of \cite{10.2307/2008251,tadmor_2003}) for a kinetic model based on \cite{Bouchut1999} and show that they yield  entropy stability for the hyperbolic system. This is in contrast to shock capturing schemes \cite{SHRINATH2023105702} based on discrete kinetic models. \\
A kinetic entropy stable scheme for continuous velocity Boltzmann's equation was recently developed in \cite{Jaiswal2022}. Although this scheme is entropy stable in the Euler limit, it employs huge number of velocities ($24^3$ for one dimensional problems) as the velocity space must be sufficiently resolved to satisfy the collision invariance. In our work, due to the usage of discrete kinetic models instead of continuous velocity Boltzmann's equation, we obtain an entropy stable scheme for the vanishing epsilon limit with very few velocities (as low as 2 for one dimensional problems). Moreover, our formalism is general enough to construct entropy stable scheme for a given hyperbolic system, while the work of \cite{Jaiswal2022} is specific to the Euler system. \\
The paper is organised as follows. In \cref{Sec: Mac model}, we briefly describe the entropy framework and entropy conservation/stability conditions required to be satisfied by a semi-discrete scheme for hyperbolic system (or macroscopic model). Then, in \cref{Sec: Vec-BGK model}, we provide a brief description of the vector-BGK model in \cite{Bouchut1999}. In \cref{Sec: Vec-kin model}, we describe our modification to vector-BGK model, termed as the vector-kinetic model. This modification allows us to obtain entropy flux potentials required for developing entropy preserving scheme for vector-kinetic model. Then, in \cref{Sec:EC,Sec:ES} we develop entropy conserving and stable schemes for vector-kinetic model, and show that these become entropy conserving and stable schemes for macroscopic model upon taking moments. In \cref{Sec:TD}, we describe the time discretisation strategies employed to complete our scheme. Then, in \cref{Sec:NR}, we verify our schemes on various numerical test problems. \Cref{Sec:Conc} concludes the paper. The list of symbols used in the paper are shown in Table \ref{tab:Symbols}. 
\begin{table}[tbhp]
\begin{center}
\begin{tabular}{|m{2cm}|m{12cm}|}
\hline
\centering \textbf{Symbol} & \centering \textbf{Description} \cr
\hline
\centering $\bold{U}$ & \centering Conserved variable vector in macroscopic model \cr
\centering $\bold{G}^{(d)}(\bold{U})$ & \centering Flux vector (along direction d) in macroscopic model \cr
\centering $\eta(\bold{U})$ & \centering Entropy function for macroscopic model \cr
\centering $\omega^{(d)}(\bold{U})$ & \centering Entropy flux function for macroscopic model  \cr
\centering $\psi^{(d)}$ & \centering Entropy flux potential for macroscopic model \cr
\centering $\bold{G}^{(d)^{\star}}_{i_{d}\pm\frac{1}{2}}$ & \centering Entropy conserving interface flux for macroscopic model \cr
\centering $\bold{Q}^{(d)^{\star}}_{i_{d}\pm\frac{1}{2}}$ & \centering Numerical viscosity corresponding to entropy conserving flux for macroscopic model \cr
\centering $\bold{G}^{(d)}_{i_{d}\pm\frac{1}{2}}$ & \centering Entropy stable interface flux for macroscopic model \cr
\centering $\bold{F}_m$ & \centering Dependent variable vector in vector-kinetic model \cr
\centering $v^{(d)}_m$ & \centering Discrete velocities in vector-kinetic model \cr
\centering $v^{(d)}_m \bold{F}_m$ & \centering Flux (along direction d) of the dependent variable vector in vector-kinetic model \cr
\centering $H_m^{\eta}$ & \centering Entropy function for vector-kinetic model \cr
\centering $v^{(d)}_m H_m^{\eta}$ & \centering Entropy flux function for vector-kinetic model \cr
\centering $\chi^{(d)}_m$ & \centering Entropy flux potential for vector-kinetic model \cr
\centering $\left(v^{(d)}_m\bold{F}_{m}\right)^{\star}_{i_{d}\pm\frac{1}{2}}$ & \centering Entropy conserving interface flux for vector-kinetic model \cr
\centering $\bold{Q}^{(d)^{\star}}_{m_{i_{d}\pm\frac{1}{2}}}$ & \centering Numerical viscosity corresponding to entropy conserving flux for vector-kinetic model \cr
\centering $\left(v^{(d)}_m\bold{F}_{m}\right)_{i_{d}\pm\frac{1}{2}}$ & \centering Entropy stable interface flux for vector-kinetic model \cr
\centering $\bold{V}$ & \centering Entropy variable \cr
\hline
\end{tabular}
\caption{\centering Table of symbols} 
\label{tab:Symbols}
\end{center}
\end{table}
\section{Macroscopic model}
\label{Sec: Mac model}
Consider the hyperbolic system (or macroscopic model),
\begin{equation}
\label{Mac model}
\partial_t \bold{U} + \partial_{x_d} \bold{G}^{(d)}(\bold{U}) = \bold{0}
\end{equation}
where $\bold{U}:\Omega \times [0,T] \rightarrow \mathbb{R}^p$ and $\bold{G}^{(d)}(\bold{U}) : \mathbb{R}^p \rightarrow \mathbb{R}^p$, with $d \in \{ 1,2,..,D\}$. Here $\Omega$ is a convex subset of $\mathbb{R}^D$. 
\subsection{Entropy framework}
Here, we briefly recall the underlying theory (presented in \cite{10.2307/2008251,tadmor_2003,TADMOR2016}) behind development of entropy conserving/stable scheme for \cref{Mac model}. \\
If the macroscopic model in \cref{Mac model} admits convex entropy-entropy flux pair $\left(\eta(\bold{U}), \omega^{(d)}(\bold{U})\right)$ that satisfies,
\begin{equation}
\label{Ent codn for mac model}
\partial_{\bold{U}} \omega^{(d)}=\partial_{\bold{U}}  \eta \cdot \partial_{\bold{U}} \bold{G}^{(d)} \Leftrightarrow \partial^2_{\bold{U}} \eta \cdot \partial_{\bold{U}}  \bold{G}^{(d)} \text{ is symmetric}
\end{equation}
then the following entropy inequality holds.
\begin{equation}
\label{Ent ineq for mac model}
\partial_t \eta(\bold{U}) + \partial_{x_d} \omega^{(d)}(\bold{U}) \leq 0
\end{equation}
Equality holds in smooth regions, while strict inequality holds in non-smooth regions. \\
Due to the convexity of $\eta(\bold{U})$,  there exists one-one correspondence $\bold{U} \rightarrow \bold{V} :=\partial_{\bold{U}}  \eta$ such that the following equivalent symmetric form of \cref{Mac model} holds true. 
\begin{equation}
\partial_{\bold{V}}\bold{U} \ \partial_t \bold{V}   + \partial_{\bold{U}}\bold{G}^{(d)} \ \partial_{\bold{V}}\bold{U} \ \partial_{x_d} \bold{V} = \bold{0}
\end{equation}
Here, $\partial_{\bold{V}}\bold{U}= \left( \partial^2_{\bold{U}} \eta(\bold{U}) \right)^{-1}$ is symmetric positive-definite (due to the convexity of $\eta(\bold{U})$) and $\partial_{\bold{V}}\bold{G}^{(d)}=\partial_{\bold{U}}\bold{G}^{(d)} \ \partial_{\bold{V}}\bold{U}$ is symmetric (refer Harten \cite{HARTEN1983151} for theorems due to Godunov and Mock). \\
Further, the compatibility condition in \cref{Ent codn for mac model} can be re-written in terms of entropy variable $\bold{V}$, thanks to the convexity of $\eta(\bold{U})$ that assures existence of $\left( \partial_{\bold{U}} \bold{V} \right)^{-1}$.
\begin{equation}
\label{Ent codn_V for mac model}
\partial_{\bold{V}}  \omega^{(d)} =\bold{V} \cdot \partial_{\bold{V}}  \bold{G}^{(d)}
\end{equation}
Due to the symmetric nature of $\partial_{\bold{V}} \bold{G}^{(d)}$,  there exist potentials $\psi^{(d)}(\bold{V})$ such that $\partial_{\bold{V}} \psi^{(d)}=\bold{G}^{(d)} (\bold{V})$. Therefore, according to \cref{Ent codn_V for mac model},  there exist entropy flux potentials,
\begin{equation}
\label{efp for mac model}
\psi^{(d)}(\bold{V})=\bold{V} \cdot \bold{G}^{(d)}(\bold{V})-\omega^{(d)}(\bold{V})
\end{equation}
\subsection{Entropy conserving scheme}
Consider a structured grid with grid size $\Delta x_d$ along each direction $d$. Then, a three-point (along each direction $d$) semi-discrete conservative scheme for \cref{Mac model} is,
\begin{equation}
\label{EC scheme for mac model}
\frac{d}{d t} \bold{U}_i + \frac{1}{\Delta x_d} \left( \bold{G}^{(d)^{\star}}_{i_{d}+\frac{1}{2}}-\bold{G}^{(d)^{\star}}_{i_{d}-\frac{1}{2}}\right)=\bold{0}
\end{equation}
Here $i$ denotes the index for cell centre of each cell/finite volume, and $i_{d}\pm\frac{1}{2}$ denote indices for right/left interfaces of cell $i$ along direction $d$. For consistency, the numerical flux $\bold{G}^{(d)^{\star}}_{i_{d}\pm \frac{1}{2}} := \bold{G}^{(d)^{\star}} _{i_{d}\pm \frac{1}{2}} \left( \bold{U}_i ,\bold{U}_{i_{d}\pm1} \right)$ is such that $\bold{G}^{(d)^{\star}}_{i_{d}\pm \frac{1}{2}} (\bold{U},\bold{U}) = \bold{G}^{(d)}(\bold{U})$, where $i_{d}\pm1$ denote indices for the cell centres of cells to the right/left of cell $i$ along direction $d$.\\ 
The scheme in \cref{EC scheme for mac model} is entropy conserving iff the interface numerical fluxes satisfy the entropy conserving condition (derived in \cite{10.2307/2008251,tadmor_2003}), 
\begin{equation}
\label{EC condn for mac model}
\left< \left[\left[\bold{V}\right]\right]_{i_{d}+\frac{1}{2}}, \bold{G}^{(d)^{\star}}_{i_{d}+\frac{1}{2}} \right>= \left[\left[\psi^{(d)}\right]\right]_{i_{d}+\frac{1}{2}}
\end{equation}
Here, $\left[\left[(.)\right]\right]_{i_{d}+\frac{1}{2}}$ denotes the jump $(.)_{i_{d}+1}-(.)_i$. Then, the following three-point (along each direction $d$) entropy equality holds true.
\begin{equation}
\label{EC eqn for mac model}
\frac{d}{d t} \eta \left( \bold{V}_i \right) + \frac{1}{\Delta x_d} \left( \omega^{(d)^{\star}}_{i_{d}+\frac{1}{2}}-\omega^{(d)^{\star}}_{i_{d}-\frac{1}{2}} \right)=0
\end{equation}
The interface numerical entropy flux consistent with  \cref{efp for mac model} is given by
\begin{equation}
\label{EC Ent flux for mac model}
\omega^{(d)^{\star}}_{i_{d}\pm\frac{1}{2}} = \frac{1}{2} \left( \bold{V}_i + \bold{V}_{i_{d}\pm1} \right) \cdot \bold{G}^{(d)^{\star}}_{i_{d}\pm\frac{1}{2}} - \frac{1}{2} \left(  \psi^{(d)}_i + \psi^{(d)}_{i_{d}\pm1}\right) 
\end{equation}
Further, the entropy conserving numerical flux $\bold{G}^{(d)^{\star}}_{i_{d}+\frac{1}{2}}$ satisfying \cref{EC condn for mac model} can be evaluated along the path $\bold{V}_{i_{d}+\frac{1}{2}} \left( \xi \right) = \bold{V}_i + \xi \Delta \bold{V}_{i_{d}+\frac{1}{2}}$ as,
\begin{equation}
\bold{G}^{(d)^{\star}}_{i_{d}\pm\frac{1}{2}} = \frac{1}{2} \left( \bold{G}^{(d)}_i + \bold{G}^{(d)}_{i_{d}\pm1} \right) -\frac{1}{2} \bold{Q}^{(d)^{\star}}_{i_{d}\pm\frac{1}{2}} \left[\left[ \bold{V} \right]\right]_{i_{d}\pm\frac{1}{2}} 
\end{equation}
with
\begin{equation}
\label{Q for mac model}
\bold{Q}^{(d)^{\star}}_{i_{d}+\frac{1}{2}} = \int_{0}^1 \left( 2 \xi -1 \right) \partial_{\bold{V}} \bold{G}^{(d)} \left( \bold{V}_{i_{d}+\frac{1}{2}} \left( \xi \right) \right) d \xi
\end{equation}
The term $\bold{Q}^{(d)^{\star}}_{i_{d}\pm\frac{1}{2}}$ which is symmetric (need not be positive-definite) is considered as numerical viscosity coefficient matrix. This counterbalances dispersion from the average flux. Further, the entropy conserving scheme is second order accurate in space (refer \cite{10.2307/2008251,tadmor_2003}). Construction of higher order entropy conserving fluxes as linear combinations of second order accurate entropy conserving fluxes $\bold{G}^{(d)^{\star}}_{i_{d}\pm\frac{1}{2}}$ is discussed in \cite{doi:10.1137/S003614290240069X}. 
\subsection{Entropy stable scheme}
The three-point (along each direction $d$) consistent flux,
\begin{equation}
\label{ES flux for mac model}
\bold{G}^{(d)}_{i_{d}\pm\frac{1}{2}} = \bold{G}^{(d)^{\star}}_{i_{d}\pm\frac{1}{2}} -\frac{1}{2} \bold{D}^{(d)}_{i_{d}\pm\frac{1}{2}} \left[\left[ \bold{V} \right]\right]_{i_{d}\pm\frac{1}{2}}
\end{equation}
with $\bold{D}^{(d)}_{i_{d}\pm\frac{1}{2}} = \bold{Q}^{(d)}_{i_{d}\pm\frac{1}{2}} - \bold{Q}^{(d)^{\star}}_{i_{d}\pm\frac{1}{2}}$ is entropy stable if and only if $\bold{D}^{(d)}_{i_{d}\pm\frac{1}{2}}$ is positive-definite. Here $\bold{Q}^{(d)}_{i_{d}\pm\frac{1}{2}}$ is the numerical viscosity coefficient matrix corresponding to entropy stable scheme. The scheme then satisfies the three-point entropy inequality,
\begin{equation}
\label{ES eqn for mac model}
\frac{d}{d t} \eta \left( \bold{V}_i \right) + \frac{1}{\Delta x_d} \left( \omega^{(d)}_{i_{d}+\frac{1}{2}}-\omega^{(d)}_{i_{d}-\frac{1}{2}} \right) = -\frac{1}{4\Delta x_d}  \left( \left[\left[ \bold{V} \right]\right]_{i_{d}+\frac{1}{2}} \cdot \bold{D}^{(d)}_{i_{d}+\frac{1}{2}} \left[\left[ \bold{V} \right]\right]_{i_{d}+\frac{1}{2}} + \left[\left[ \bold{V} \right]\right]_{i_{d}-\frac{1}{2}} \cdot \bold{D}^{(d)}_{i_{d}-\frac{1}{2}} \left[\left[ \bold{V} \right]\right]_{i_{d}-\frac{1}{2}} \right) \leq 0   
\end{equation}
Here, the consistent numerical entropy flux at interface is given by,
\begin{equation}
\label{ES Ent flux for mac model}
\omega^{(d)}_{i_{d}+\frac{1}{2}} = \omega^{(d)^{\star}}_{i_{d}+\frac{1}{2}} - \frac{1}{4} \left( \bold{V}_i + \bold{V}_{i_{d}+1}\right) \cdot  \bold{D}^{(d)}_{i_{d}+\frac{1}{2}}\left[\left[ \bold{V}\right]\right]_{i_{d}+\frac{1}{2}}
\end{equation}
The entropy stable flux $\bold{G}^{(d)}_{i_{d}\pm\frac{1}{2}}$ given by \cref{ES flux for mac model} is first order accurate in space (refer Tadmor \cite{10.2307/2008251,tadmor_2003}). To achieve higher order accuracy in space, the term $\left[\left[ \bold{V} \right]\right]_{i_{d}+\frac{1}{2}}$ in \cref{ES flux for mac model} must be replaced by $\left<\left< \bold{V} \right>\right>_{i_{d}+\frac{1}{2}}=\bold{V}^-_{i_{d}+1}-\bold{V}^+_i$ where $\bold{V}^-_{i_{d}+1}$ and $\bold{V}^+_i$ are higher order reconstructions of $\bold{V}$ at interface $i_{d}+\frac{1}{2}$ (refer \cite{doi:10.1137/110836961}).

\section{Vector-BGK model}
\label{Sec: Vec-BGK model}
In this section, we briefly describe the vector-BGK model presented in \cite{Bouchut1999}. Consider,
\begin{equation}
\label{vec bgk model}
\partial_t \bold{f}_m + \partial_{x_d} \left( v^{(d)}_m\bold{f}_m \right) = - \frac{1}{\epsilon} \left( \bold{f}_m  - \bold{F}_m(\bold{U}) \right)
\end{equation}
where $\epsilon$ is the relaxation parameter. Here, $\bold{f}_m:= \bold{f}_m(x_1,..,x_d,..,x_D,v^{(1)}_m,..,v^{(d)}_m,..,v^{(D)}_m,t) \in \mathbb{R}^p$, $\bold{F}_m:\mathbb{R}^p \rightarrow \mathbb{R}^p$, $m \in \{1,2,..,M \}$ and $M$ is the number of discrete velocities. Splitting of streaming and relaxation operators in \cref{vec bgk model} gives,
\begin{eqnarray}
\text{Streaming:} & \partial_t \bold{f}_m + \partial_{x_d} \left( v^{(d)}_m\bold{f}_m \right) = \bold{0} \\ 
\text{Relaxation:} & \frac{d}{dt} \bold{f}_m = - \frac{1}{\epsilon} \left( \bold{f}_m - \bold{F}_m(\bold{U}) \right)
\end{eqnarray}
Instantaneous relaxation (\text{i.e.,} $\epsilon=0$) in the relaxation equation above  yields $\bold{f}_m=\bold{F}_m(\bold{U})$. This is inserted into the streaming equation for its evolution. Now, it can be seen that if the following relations are satisfied,
\begin{equation}
\label{mom vec kin model}
\sum_{m=1}^M \bold{F}_m(\bold{U}) = \bold{U} \text{ and } \sum_{m=1}^M v^{(d)}_m \bold{F}_m(\bold{U}) = \bold{G}^{(d)}(\bold{U})
\end{equation} 
then $\sum_{m=1}^M $ \cref{vec bgk model} $\rightarrow$ \cref{Mac model} as $\epsilon \rightarrow 0$. \\
\subsection{Entropy framework}
Following the definition of entropy function for vector-BGK model given by equations (E0)-(E2) in \cite{Bouchut1999}, let us define the entropy function $H_m^{\eta}\left( \bold{f}_m\right)$ as:
\begin{eqnarray}
\label{def_H1} H_m^{\eta}\left( \bold{f}_m\right) \text{ is a convex function with respect to } \bold{f}_m \\
\label{def_H2} \sum_{m=1}^M H_m^{\eta}\left( \bold{F}_m(\bold{U}) \right) = \eta(\bold{U}) \\
\label{def_H3} \sum_{m=1}^M H_m^{\eta}\left( \bold{F}_m(\bold{U}) \right) \leq \sum_{m=1}^M H_m^{\eta}\left( \bold{f}_m \right) 
\end{eqnarray}
Then, taking inner product of \cref{vec bgk model} with the sub-differential of $H_m^{\eta}$ at $\bold{F}_m(\bold{U})$ and using (\ref{def_H1}), (\ref{def_H2}) and (\ref{def_H3}), the following is obtained. 
\begin{multline}
\label{H eqn vec bgk model}
\partial_t H_m^{\eta} (\bold{f}_m) + \partial_{x_d} \left( v^{(d)}_m H_m^{\eta} (\bold{f}_m)  \right) \leq \frac{1}{\epsilon} \biggl( H_m^{\eta} \left( \bold{F}_m(\bold{U}) \right) - H_m^{\eta} \left( \bold{f}_m \right) \biggr) \\
\Rightarrow \sum_{m=1}^M \biggl( \partial_t H_m^{\eta} (\bold{f}_m) + \partial_{x_d} \left( v^{(d)}_m H_m^{\eta} (\bold{f}_m)\right)  \biggr) \leq 0 \\
\Rightarrow \partial_t \eta(\bold{U}) + \partial_{x_d}  \left( \sum_{m=1}^M  v^{(d)}_m H_m^{\eta}  (\bold{F}_m(\bold{U})) \right)\leq 0 \text{ in the limit } \epsilon \rightarrow 0
\end{multline}
If $\omega^{(d)}(\bold{U})=\sum_{m=1}^M  v^{(d)}_m H_m^{\eta}  (\bold{F}_m(\bold{U}))$, then \cref{H eqn vec bgk model} is same as \cref{Ent ineq for mac model}. The reader is referred to \cite{Bouchut1999} for details.\\
Thus, entropy inequality of the macroscopic model (\cref{Mac model}) can be obtained as minimisation of entropies of the vector-BGK model (\cref{vec bgk model}). This inspires one to develop entropy structure preserving numerical schemes for vector-BGK model that recover the entropy inequality of equivalent macroscopic scheme. However, the framework of vector-BGK model does not ensure the existence of $ \partial^2_{\bold{f}_m} H_m^{\eta} \left( \bold{F}_m(\bold{U}) \right)$ which is crucial in obtaining entropy flux potentials that allow for the consistent definition of interface numerical entropy fluxes. Hence, we resort to a much simpler model in the relaxed limit without the stiff relaxation parameter (hereafter referred as \textit{vector-kinetic model}), and make the necessary modification to allow for the definition of entropy flux potentials. 

\section{Vector-kinetic model}
\label{Sec: Vec-kin model}
In this model, we consider the evolution of relaxed limit ($\epsilon=0$):
\begin{equation}
\label{vec kin model}
\partial_t \bold{F}_m + \partial_{x_d} \left( v^{(d)}_m\bold{F}_m \right) = \bold{0}
\end{equation}
Let us define $\bold{F}_m(\bold{U})$ as in \cite{Bouchut1999},
\begin{equation}
\label{F def vec kin model}
\bold{F}_m (\bold{U})=a_m\bold{U}+b^{(d)}_m\bold{G}^{(d)}\left( \bold{U} \right)
\end{equation}
with 
\begin{gather}
\label{Mom cons. vec kin model_1}
\sum_{m=1}^M a_m = 1, \ \sum_{m=1}^M b^{(d)}_m = 0 \\  \sum_{m=1}^M v^{(j)}_m a_m  = 0, \ \sum_{m=1}^M v^{(j)}_m b^{(d)}_m  = \delta_{jd} \label{Mom cons. vec kin model_2}
\end{gather}
In the light of moment constraints in \cref{Mom cons. vec kin model_1,Mom cons. vec kin model_2}, the definition of $\bold{F}_m(\bold{U})$ in \cref{F def vec kin model} satisfies \cref{mom vec kin model}. 
\subsection{Entropy framework}
Define $H_m^{\eta}$ as in \cite{Bouchut1999},
\begin{equation}
\label{H def vec kin model}
H_m^{\eta} (\bold{U}) = a_m \eta \left( \bold{U}\right) + b^{(d)}_m \omega^{(d)} \left( \bold{U}\right)
\end{equation}
Due to the constraints in \cref{Mom cons. vec kin model_1,Mom cons. vec kin model_2}, $H_m^{\eta}$ satisfies,
\begin{equation}
\sum_{m=1}^M H_m^{\eta} (\bold{U}) = \eta(\bold{U}) \text{ and } \sum_{m=1}^M v^{(d)}_m H_m^{\eta}(\bold{U}) = \omega^{(d)}(\bold{U})
\end{equation}
We assume that the eigenvalues of $\partial_{\bold{U}} \bold{F}_m$ are positive, unlike in \cite{Bouchut1999} where the eigenvalues are considered to be non-negative. It will be seen that this modification allows the definition of entropy flux potentials required in the construction of entropy preserving numerical scheme. As $\partial_{\bold{U}} \bold{F}_m$ is now invertible, $\partial_{\bold{F}_m} H_m^{\eta}$ satisfying $\partial_{\bold{U}} H^{\eta}_m = \partial_{\bold{F}_m} H^{\eta}_m \cdot  \partial_{\bold{U}} \bold{F}_m$ exists. Therefore, the inner product of \cref{vec kin model} with $\partial_{\bold{F}_m} H_m^{\eta}$ gives,
\begin{equation}
\label{Ent eq for vec kin model}
\partial_t H^{\eta}_m + \partial_{x_d} \left( v^{(d)}_m H^{\eta}_m \right) = 0
\end{equation}
It can be seen that $\sum_{m=1}^M \left( \cref{Ent eq for vec kin model}\right)$ becomes \cref{Ent ineq for mac model} with equality. Motivated by this, in this paper, we develop entropy preserving scheme for vector-kinetic model that recovers entropy preservation of equivalent macroscopic scheme.  
\begin{lemma}
\label{Exis ent var for vec kin model}
If $\bold{F}_m (\bold{U})$ and $H_m^{\eta} (\bold{U})$ respectively follow \cref{F def vec kin model,H def vec kin model} with constants $a_m,\ b^{(d)}_m$ satisfying the moment constraints in \cref{Mom cons. vec kin model_1,Mom cons. vec kin model_2} and rendering the eigenvalues of $\partial_{\bold{U}} \bold{F}_m$ to be positive, then $\partial_{\bold{F}_m} H_m^{\eta}=\partial_{\bold{U}}  \eta$.
\end{lemma}
\begin{proof}
Due to the compatibility condition in \cref{Ent codn for mac model}, it can be seen from differentiation (with respect to $\bold{U}$) of \cref{F def vec kin model,H def vec kin model} that $\partial_{\bold{U}} H^{\eta}_m = \partial_{\bold{U}} \eta \cdot \partial_{\bold{U}} \bold{F}_m$. Since $\partial_{\bold{U}} \bold{F}_m$ is invertible, $\partial_{\bold{U}} \eta=\partial_{\bold{U}} H^{\eta}_m \cdot \left(\partial_{\bold{U}} \bold{F}_m\right)^{-1}$. We already saw that $\partial_{\bold{F}_m} H_m^{\eta}=\partial_{\bold{U}} H^{\eta}_m \cdot \left(\partial_{\bold{U}} \bold{F}_m\right)^{-1}$.
\end{proof}
This lemma shows that the entropy variables for macroscopic and vector-kinetic models are equal, \textit{i.e.,} 
\begin{equation}
\label{equality of ent var}
    \bold{V}=\partial_{\bold{U}} \eta=\partial_{\bold{F}_m} H_m^{\eta}.
\end{equation}
The choice of constants $a_m,\ b^{(d)}_m$ satisfying assumptions in the above lemma are discussed in \ref{app constants a,b}. \\
As a consequence of \cref{Exis ent var for vec kin model}, we have $\partial^2_{\bold{F}_m} H^{\eta}_m =  \partial^2_{\bold{U}} \eta \cdot \left( \partial_{\bold{U}} \bold{F}_m \right)^{-1}$. Further, $\left( \partial^2_{\bold{U}} \eta \right)^{-1} \partial^2_{\bold{F}_m} H^{\eta}_m =   \left( \partial_{\bold{U}} \bold{F}_m \right)^{-1}$ can be expressed as 
\begin{equation}
\label{sim eqn for vec kin model}
\left( \partial^2_{\bold{U}} \eta \right)^{-\frac{1}{2}}  \left( \partial^2_{\bold{U}} \eta \right)^{-\frac{1}{2}} \left( \partial^2_{\bold{F}_m} H^{\eta}_m \right) \left( \partial^2_{\bold{U}} \eta \right)^{-\frac{1}{2}} \left( \partial^2_{\bold{U}} \eta \right)^{\frac{1}{2}} = \left( \partial_{\bold{U}} \bold{F}_m \right)^{-1}
\end{equation}
thanks to the positive-definiteness of $\partial^2_{\bold{U}} \eta$. Thus, $\left( \partial^2_{\bold{U}} \eta \right)^{-\frac{1}{2}} \left( \partial^2_{\bold{F}_m} H^{\eta}_m \right) \left( \partial^2_{\bold{U}} \eta \right)^{-\frac{1}{2}}$ and $\left( \partial_{\bold{U}} \bold{F}_m \right)^{-1}$ are similar and therefore their eigenvalues are same. 
\begin{lemma}
\label{pos def lem for vec kin model}
If $\partial^2_{\bold{U}} \eta$ is positive-definite and \cref{sim eqn for vec kin model} holds true, then $\partial^2_{\bold{F}_m} H^{\eta}_m$ is positive-definite iff the eigenvalues of $\left( \partial_{\bold{U}} \bold{F}_m \right)^{-1}$ are positive.
\end{lemma}
\begin{proof}
$\left( \partial^2_{\bold{U}} \eta \right)^{-\frac{1}{2}} \left( \partial^2_{\bold{F}_m} H^{\eta}_m \right) \left( \partial^2_{\bold{U}} \eta \right)^{-\frac{1}{2}}$ is symmetric as $\partial^2_{\bold{U}} \eta $ and $\partial^2_{\bold{F}_m} H^{\eta}_m$ are symmetric. Further, we have $\forall \bold{y}\neq\bold{0} \in \mathbb{R}^p$,
\begin{equation}
\label{pos def for vec kin model}
\bold{y} \cdot \left( \partial^2_{\bold{U}} \eta \right)^{-\frac{1}{2}} \left( \partial^2_{\bold{F}_m} H^{\eta}_m \right) \left( \partial^2_{\bold{U}} \eta \right)^{-\frac{1}{2}} \bold{y} = \bold{z} \cdot \left( \partial^2_{\bold{F}_m} H^{\eta}_m \right)\bold{z}
\end{equation}
where $\bold{z}=\left( \partial^2_{\bold{U}} \eta \right)^{-\frac{1}{2}} \bold{y} \neq \bold{0}$ (as $\partial^2_{\bold{U}} \eta$ is positive-definite). \\
$\Leftarrow$ If the eigenvalues of $\left( \partial_{\bold{U}} \bold{F}_m \right)^{-1}$ are positive, then $\left( \partial^2_{\bold{U}} \eta \right)^{-\frac{1}{2}} \left( \partial^2_{\bold{F}_m} H^{\eta}_m \right) \left( \partial^2_{\bold{U}} \eta \right)^{-\frac{1}{2}}$ is positive-definite due to \cref{sim eqn for vec kin model}. Then $\partial^2_{\bold{F}_m} H^{\eta}_m$ is rendered positive-definite by \cref{pos def for vec kin model}. \\
$\Rightarrow$ If $\partial^2_{\bold{F}_m} H^{\eta}_m$ is positive-definite, then by \cref{pos def for vec kin model} $\left( \partial^2_{\bold{U}} \eta \right)^{-\frac{1}{2}} \left( \partial^2_{\bold{F}_m} H^{\eta}_m \right) \left( \partial^2_{\bold{U}} \eta \right)^{-\frac{1}{2}}$ is positive-definite. Then, the eigenvalues of $\left( \partial_{\bold{U}} \bold{F}_m \right)^{-1}$ are positive due to \cref{sim eqn for vec kin model}. 
\end{proof}
Thus, as consequence of \cref{Exis ent var for vec kin model} and \cref{pos def lem for vec kin model}, \cref{equality of ent var} and positive-definiteness of $\partial^2_{\bold{F}_m} H^{\eta}_m$ are guaranteed iff the eigenvalues of $\partial_{\bold{U}} \bold{F}_m $ are positive. Using the one-to-one correspondence between $\bold{U}$ and $\bold{V}$, we consider $\bold{F}_m \left( \bold{U}\right) = \bold{F}_m \left( \bold{U} \left( \bold{V} \right) \right)$. Hence the vector-kinetic model in \cref{vec kin model} can be expressed in the equivalent symmetric form
\begin{equation}
\partial_{\bold{V}} \bold{F}_m \partial_t \bold{V} + \partial_{\bold{V}} \left( v^{(d)}_m\bold{F}_m \right) \partial_{x_d} \bold{V} = \bold{0}
\end{equation}
Here $\partial_{\bold{V}} \bold{F}_m =  \left( \partial^2_{\bold{F}_m} H^{\eta}_m \right)^{-1}$ is symmetric positive-definite. Due to the linearity of vector-kinetic model, $\partial_{\bold{V}} \left( v^{(d)}_m\bold{F}_m \right) = v^{(d)}_m \partial_{\bold{V}} \bold{F}_m$ is symmetric. As a result,  there exist potentials $\chi^{(d)}_m(\bold{V})$ such that 
\begin{equation}
\label{pot1 for vec kin model}
\partial_{\bold{V}} \chi^{(d)}_m=  v^{(d)}_m\bold{F}_m 
\end{equation}
Further, the compatibility condition 
\begin{equation}
\label{ent cndn for vec kin model_F}
\partial_{\bold{F}_m} \left(v^{(d)}_m H^{\eta}_m\right) = \partial_{\bold{F}_m} H^{\eta}_m \cdot  \partial_{\bold{F}_m} \left(v^{(d)}_m \bold{F}_m\right) 
\end{equation}
is also satisfied rendering $H^{\eta}_m$ as the convex entropy function for vector-kinetic model. Note that this compatibility condition is always true for any convex $H^{\eta}_m$ satisfying \cref{H def vec kin model} due to the linear nature of vector-kinetic model, unlike the compatibility condition $\left(\cref{Ent codn for mac model}\right)$ for macroscopic model. In terms of $\bold{V}$, the above compatibility condition for vector-kinetic model becomes,
\begin{equation}
\label{Ent condn for vec kin model}
\partial_{\bold{V}} \left( v^{(d)}_m H^{\eta}_m \right)  = \bold{V} \cdot \partial_{\bold{V}} \left( v^{(d)}_m\bold{F}_m \right)  
\end{equation}
thanks to the inverse of $\partial_{\bold{F}_m} \bold{V}$. Therefore, due to \cref{pot1 for vec kin model,Ent condn for vec kin model},  there exist entropy flux potentials
\begin{equation}
\label{efp for vec kin model}
\chi^{(d)}_m (\bold{V}) = \bold{V} \cdot v^{(d)}_m\bold{F}_m - v^{(d)}_m H^{\eta}_m  = \partial_{\bold{F}_m} H^{\eta}_m \cdot v^{(d)}_m\bold{F}_m - v^{(d)}_m H^{\eta}_m
\end{equation}
Thus, we have obtained the entropy flux potentials that are crucial in the construction of entropy preserving numerical scheme for vector-kinetic model. 

\section{Entropy conserving scheme for vector-kinetic model}
\label{Sec:EC}
The three-point (along each direction $d$) semi-discrete conservative scheme for vector-kinetic model in \cref{vec kin model} on a structured grid is given by, 
\begin{equation}
\label{EC scheme for vec kin model}
\frac{d}{d t} \bold{F}_{m_i} + \frac{1}{\Delta x_d} \left( \left(v^{(d)}_m\bold{F}_{m}\right)^{\star}_{i_{d}+\frac{1}{2}}-\left(v^{(d)}_m\bold{F}_{m}\right)^{\star}_{i_{d}-\frac{1}{2}}\right)=\bold{0}
\end{equation}
Here, $\bold{F}_{m_i}(t)=\bold{F}_m\left(\bold{V}_i(t)\right)$ and consistent $\left(v^{(d)}_m\bold{F}_{m}\right)^{\star}_{i_{d}+\frac{1}{2}} = v^{(d)}_m\bold{F}_{m} (\bold{V}_i,\bold{V}_{i_{d}+1})$ is such that $v^{(d)}_m\bold{F}_{m} (\bold{V},\bold{V})=v^{(d)}_m\bold{F}_{m} (\bold{V})$. Consider the inner product $\left(\partial_{\bold{F}_m} H^{\eta}_m\right)_i \cdot \left(v^{(d)}_m\bold{F}_{m}\right)^{\star}_{i_{d}\pm\frac{1}{2}}$:
\begin{eqnarray*}
    \left(\partial_{\bold{F}_m} H^{\eta}_m\right)_i \cdot \left(v^{(d)}_m\bold{F}_{m}\right)^{\star}_{i_{d}\pm\frac{1}{2}} &=& \frac{1}{2} \left( \left(\partial_{\bold{F}_m} H^{\eta}_m\right)_{i_{d}\pm1} + \left(\partial_{\bold{F}_m} H^{\eta}_m\right)_i \right) \cdot \left(v^{(d)}_m\bold{F}_{m}\right)^{\star}_{i_{d}\pm\frac{1}{2}} \\ & & - \frac{1}{2} \left( \left(\partial_{\bold{F}_m} H^{\eta}_m\right)_{i_{d}\pm1} - \left(\partial_{\bold{F}_m} H^{\eta}_m\right)_i \right) \cdot \left(v^{(d)}_m\bold{F}_{m}\right)^{\star}_{i_{d}\pm\frac{1}{2}}  
\end{eqnarray*}
If the interface numerical flux $\left(v^{(d)}_m\bold{F}_{m}\right)^{\star}_{i_{d}+\frac{1}{2}}$ satisfies the entropy conserving condition,
\begin{equation}
\label{EC condn for vec kin model}
\left< \left[\left[\partial_{\bold{F}_m} H^{\eta}_m\right]\right]_{i_{d}+\frac{1}{2}}, \left(v^{(d)}_m\bold{F}_{m}\right)^{\star}_{i_{d}+\frac{1}{2}} \right>= \left[\left[\chi^{(d)}_{m}\right]\right]_{i_{d}+\frac{1}{2}}
\end{equation}
then,
\begin{eqnarray*}
    \left(\partial_{\bold{F}_m} H^{\eta}_m\right)_i \cdot \left(v^{(d)}_m\bold{F}_{m}\right)^{\star}_{i_{d}\pm\frac{1}{2}} &=& \frac{1}{2} \left( \left(\partial_{\bold{F}_m} H^{\eta}_m\right)_{i_{d}\pm1} + \left(\partial_{\bold{F}_m} H^{\eta}_m\right)_i \right) \cdot \left(v^{(d)}_m\bold{F}_{m}\right)^{\star}_{i_{d}\pm\frac{1}{2}} - \frac{1}{2} \left( \chi^{(d)}_{m_{i_{d}\pm1}} - \chi^{(d)}_{m_i} \right)
\end{eqnarray*}
Thus, the inner product of \cref{EC scheme for vec kin model} with $\left(\partial_{\bold{F}_m} H^{\eta}_m\right)_i$ gives the three-point entropy equality, 
\begin{equation}
\label{EC eqn for vec kin model}
\frac{d}{d t} H^{\eta}_{m_i}  + \frac{1}{\Delta x_d} \left( \left(v^{(d)}_mH^{\eta}_{m}\right)^{\star}_{i_{d}+\frac{1}{2}}-\left(v^{(d)}_mH^{\eta}_{m}\right)^{\star}_{i_{d}-\frac{1}{2}}\right)=0
\end{equation}
iff it satisfies \cref{EC condn for vec kin model}, and the interface numerical entropy fluxes $\left(v^{(d)}_mH^{\eta}_{m}\right)^{\star}_{i_{d}\pm\frac{1}{2}}$ consistent with \cref{efp for vec kin model} are given by,
\begin{equation}
\label{EC ef for vec kin model}
\left(v^{(d)}_mH^{\eta}_{m}\right)^{\star}_{i_{d}\pm\frac{1}{2}} = \frac{1}{2} \left( \left( \partial_{\bold{F}_m} H^{\eta}_m\right)_i + \left( \partial_{\bold{F}_m} H^{\eta}_m\right)_{i_{d}\pm1} \right) \cdot \left(v^{(d)}_m\bold{F}_{m}\right)^{\star}_{i_{d}\pm\frac{1}{2}} - \frac{1}{2} \left(  \chi^{(d)}_{m_i} + \chi^{(d)}_{m_{i_{d}\pm1}}\right)
\end{equation}
It is seen that the entropy flux potentials $\chi^{(d)}_{m_i}$ enable us to consistently relate the two interfacial unknowns, numerical fluxes $\left(v^{(d)}_m\bold{F}_{m}\right)^{\star}_{i_{d}\pm\frac{1}{2}}$ and numerical entropy fluxes $\left(v^{(d)}_mH^{\eta}_{m}\right)^{\star}_{i_{d}\pm\frac{1}{2}}$. 
Further, let us define the interface numerical fluxes for macroscopic model as the moment of interface numerical fluxes for vector-kinetic model as,
\begin{equation}
\label{int flux mac model- vec kin model}
\bold{G}^{(d)^{\star}}_{i_{d}\pm\frac{1}{2}} = \sum_{m=1}^M \left(v^{(d)}_m\bold{F}_{m}\right)^{\star}_{i_{d}\pm\frac{1}{2}}
\end{equation}
\begin{theorem}
\label{thm EC for vec kin-mac model}
If the three-point semi-discrete conservative scheme $(\cref{EC scheme for vec kin model})$ for vector-kinetic model with 
\begin{itemize}
\item{$\bold{F}_{m_i}=a_m\bold{U}_i+b^{(d)}_m \bold{G}^{(d)}_i, \ \forall i$} \item{interface numerical fluxes $\left(v^{(d)}_m\bold{F}_{m}\right)^{\star}_{i_{d}\pm\frac{1}{2}}$ satisfying the entropy conserving condition in \cref{EC condn for vec kin model} and}
\item{constants $a_m$, $b^{(d)}_m$ satisfying the moment constraints in \cref{Mom cons. vec kin model_1,Mom cons. vec kin model_2} while rendering positivity of eigenvalues of $\partial_{\bold{U}} \bold{F}_m$}
\end{itemize}
is used, and if the convex entropy function corresponding to it is $H^{\eta}_{m_i}=a_m\eta_i+b^{(d)}_m \omega^{(d)}_i, \ \forall i$, then 
\begin{enumerate}
\item{$\sum_{m=1}^M \left(\cref{EC scheme for vec kin model}\right)$ becomes 
\begin{equation}
\label{thm EC 1st item}
\frac{d}{d t} \bold{U}_i + \frac{1}{\Delta x_d} \left( \bold{G}^{(d)^{\star}}_{i_{d}+\frac{1}{2}}-\bold{G}^{(d)^{\star}}_{i_{d}-\frac{1}{2}}\right)=\bold{0}
\end{equation} 
with $\bold{G}^{(d)^{\star}}_{i_{d}\pm\frac{1}{2}}$ given by \cref{int flux mac model- vec kin model},}
\item{the interface numerical flux $\bold{G}^{(d)^{\star}}_{i_{d}\pm\frac{1}{2}}$ given by \cref{int flux mac model- vec kin model} satisfies the entropy conserving condition for macroscopic model $\left(\cref{EC condn for mac model}\right)$, and}
\item{the three-point entropy equality for macroscopic model $\left(\cref{EC eqn for mac model}\right)$ holds true with interface numerical entropy flux $\omega^{(d)^{\star}}_{i_{d}\pm\frac{1}{2}}$ given by \cref{EC Ent flux for mac model}.}
\end{enumerate}
\end{theorem}
\begin{proof}
Due to moment constraint in \cref{Mom cons. vec kin model_1}, $\sum_{m=1}^M \bold{F}_{m_i}=\bold{U}_i$. Therefore, $\sum_{m=1}^M \left( \cref{EC scheme for vec kin model}\right)$ becomes \cref{thm EC 1st item} with $\bold{G}^{(d)^{\star}}_{i_{d}\pm\frac{1}{2}}$ given by \cref{int flux mac model- vec kin model}, thus proving 1.  \\
By \cref{equality of ent var}, $\left[\left[\partial_{\bold{F}_m} H^{\eta}_m \right]\right]_{i_{d}\pm\frac{1}{2}}=\left[\left[\bold{V}\right]\right]_{i_{d}\pm\frac{1}{2}}=  \left[\left[\partial_{\bold{U}}  \eta\right]\right]_{i_{d}\pm\frac{1}{2}}$ is not a function of $m$. Hence, the moment of \cref{EC condn for vec kin model} gives,
\begin{equation}
\label{ent condn vec kin-mac model}
\left< \left[\left[\bold{V}\right]\right]_{i_{d}\pm\frac{1}{2}} , \sum_{m=1}^M \left( v^{(d)}_m \bold{F}_m \right)^{\star}_{i_{d}\pm\frac{1}{2}} \right> = \left[\left[\sum_{m=1}^M \chi^{(d)}_m\right]\right]_{i_{d}\pm\frac{1}{2}}
\end{equation}
From \cref{efp for vec kin model}, it can be seen that $\chi^{(d)}_{m_i} = \bold{V}_i .v^{(d)}_m\bold{F}_{m_i} - v^{(d)}_m H^{\eta}_{m_i},\ \forall i$. Hence, $\sum_{m=1}^M \chi^{(d)}_{m_i} = \bold{V}_i . \sum_{m=1}^M \left( v^{(d)}_m\bold{F}_{m_i} \right) - \sum_{m=1}^M \left( v^{(d)}_m H^{\eta}_{m_i} \right), \ \forall i $. We also have $\sum_{m=1}^M v^{(d)}_m \bold{F}_{m_i}=\bold{G}^{(d)}_i$ and $\sum_{m=1}^M v^{(d)}_m H^{\eta}_{m_i}= \omega^{(d)}_i, \ \forall i$ due to the action of moment constraint in \cref{Mom cons. vec kin model_2} on $\bold{F}_{m_i}$ and $H^{\eta}_{m_i}$. Therefore, by \cref{efp for mac model}, $\sum_{m=1}^M \chi^{(d)}_{m_i}=\psi^{(d)}_i, \ \forall i$. Using this and \cref{int flux mac model- vec kin model} in \cref{ent condn vec kin-mac model}, we obtain,
\begin{equation}
\left< \left[\left[\bold{V}\right]\right]_{i_{d}\pm\frac{1}{2}} , \bold{G}^{(d)^{\star}}_{i_{d}\pm\frac{1}{2}} \right> = \left[\left[\psi^{(d)}\right]\right]_{i_{d}\pm\frac{1}{2}}
\end{equation}
This proves 2. \\
We know that the three-point entropy equality in \cref{EC eqn for vec kin model} holds true corresponding to the assumptions stated in \cref{thm EC for vec kin-mac model}. Since $\sum_{m=1}^M (H^{\eta}_m)_i=\eta_i, \ \forall i$ (due to the action of moment constraint in \cref{Mom cons. vec kin model_1} on $(H^{\eta}_m)_i$), moment of \cref{EC eqn for vec kin model} gives,
\begin{equation}
\frac{d}{d t} \eta_i  + \frac{1}{\Delta x_d} \left(  \sum_{m=1}^M \left(v^{(d)}_mH^{\eta}_{m}\right)^{\star}_{i_{d}+\frac{1}{2}}-  \sum_{m=1}^M \left(v^{(d)}_mH^{\eta}_{m}\right)^{\star}_{i_{d}-\frac{1}{2}}\right)=0 
\end{equation}
Since $\left(\partial_{\bold{F}_m} H^{\eta}_m\right)_i=\bold{V}_i=\left(\partial_{\bold{U}}  \eta\right)_i$ is not a function of $m$ (by \cref{equality of ent var}), moment of $\left(v^{(d)}_m H^{\eta}_{m}\right)^{\star}_{i_{d}\pm\frac{1}{2}}$ given by \cref{EC ef for vec kin model} yields, 
\begin{equation}
\sum_{m=1}^M \left(v^{(d)}_mH^{\eta}_{m}\right)^{\star}_{i_{d}\pm\frac{1}{2}} = \frac{1}{2} \left( \bold{V}_i + \bold{V}_{i_{d}\pm1} \right) \cdot \sum_{m=1}^M  \left(v^{(d)}_m\bold{F}_{m}\right)^{\star}_{i_{d}\pm\frac{1}{2}}  - \frac{1}{2} \left(  \sum_{m=1}^M \chi^{(d)}_{m_i} + \sum_{m=1}^M \chi^{(d)}_{m_{i_{d}\pm1}}\right)
\end{equation}
We have already seen that $\sum_{m=1}^M \chi^{(d)}_{m_i}=\psi^{(d)}_i, \ \forall i$. Using this and \cref{int flux mac model- vec kin model}, we obtain,
\begin{equation}
\sum_{m=1}^M \left(v^{(d)}_mH^{\eta}_{m}\right)^{\star}_{i_{d}\pm\frac{1}{2}} = \frac{1}{2} \left( \bold{V}_i + \bold{V}_{i_{d}\pm1} \right) \cdot \bold{G}^{(d)^{\star}}_{i_{d}\pm\frac{1}{2}} - \frac{1}{2} \left(   \psi^{(d)}_i +  \psi^{(d)}_{i_{d}\pm1}\right)
\end{equation}
It can be seen from \cref{EC Ent flux for mac model} that $\sum_{m=1}^M \left(v^{(d)}_mH^{\eta}_{m}\right)^{\star}_{i_{d}\pm\frac{1}{2}}= \omega^{(d)^{\star}}_{i_{d}\pm\frac{1}{2}}$. This proves 3. 
\end{proof}
In the light of \cref{equality of ent var} resulting from \cref{Exis ent var for vec kin model}, moments involved in the proof of above theorem become linear since $\partial_{\bold{F}_m} H_m^{\eta}$ is not a function of $m$. This plays a pivotal role in showing that entropy conserving scheme for vector-kinetic model results in an entropy conserving scheme for macroscopic model. 
\begin{remark}
In the above proof, the three-point entropy equality for macroscopic model $\left(\cref{EC eqn for mac model}\right)$ with interface numerical entropy flux $\omega^{(d)^{\star}}_{i_{d}\pm\frac{1}{2}}$ given by \cref{EC Ent flux for mac model} is obtained as moment of three-point entropy equality for vector-kinetic model. Unlike this, we can also obtain \cref{EC eqn for mac model} directly at the macroscopic level as a consequence of $\bold{G}^{(d)^{\star}}_{i_{d}\pm\frac{1}{2}}= \sum_{m=1}^M \left(v^{(d)}_m\bold{F}_{m}\right)^{\star}_{i_{d}\pm\frac{1}{2}}$ satisfying the entropy conserving condition for macroscopic model $\left(\cref{EC condn for mac model}\right)$.  
\end{remark}
The entropy conserving fluxes satisfying \cref{EC condn for vec kin model} can be evaluated using an integral along the path $\bold{V}_{i_{d}+\frac{1}{2}} \left( \xi \right) = \bold{V}_i + \xi \Delta \bold{V}_{i_{d}+\frac{1}{2}}$ as,
\begin{equation}
\label{avg+numdiff form for vec kin model}
\left(v^{(d)}_m\bold{F}_{m}\right)^{\star}_{i_{d}\pm\frac{1}{2}} = \int_{0}^1 \left(v^{(d)}_m\bold{F}_{m}\right) \left( \bold{V}_{i_{d}+\frac{1}{2}} \left( \xi \right) \right) d\xi  = \frac{1}{2} \left( v^{(d)}_m\bold{F}_{m_i} + v^{(d)}_m\bold{F}_{m_{i_{d}\pm1}} \right) -\frac{1}{2} \bold{Q}^{(d)^{\star}}_{m_{i_{d}\pm\frac{1}{2}}} \left[\left[ \bold{V} \right]\right]_{i_{d}\pm\frac{1}{2}}
\end{equation}
where 
\begin{equation}
\label{num visc of ec_kin}
\bold{Q}^{(d)^{\star}}_{m_{i_{d}+\frac{1}{2}}} = \int_{0}^1 \left( 2 \xi -1 \right) \partial_{\bold{V}} \left(v^{(d)}_m\bold{F}_{m}\right) \left( \bold{V}_{i_{d}+\frac{1}{2}} \left( \xi \right) \right) d \xi
\end{equation}
Although $\partial_{\bold{V}} \left(v^{(d)}_m\bold{F}_{m}\right) \left( \bold{V}_{i_{d}+\frac{1}{2}} \left( \xi \right) \right)$ is symmetric positive-definite, the term $\bold{Q}^{(d)^{\star}}_{m_{i_{d}+\frac{1}{2}}}$ is only symmetric (need not be positive-definite). This is considered as numerical viscosity coefficient matrix that counterbalances the dispersion from average flux. Integration by parts of $\bold{Q}^{(d)^{\star}}_{m_{i_{d}+\frac{1}{2}}}$ yields,
\begin{equation} 
\bold{Q}^{(d)^{\star}}_{m_{i_{d}+\frac{1}{2}}}= \int_{0}^1 \left( \xi -\xi^2 \right) \partial_{\bold{VV}} \left(v^{(d)}_m\bold{F}_{m}\right)\left( \bold{V}_{i_{d}+\frac{1}{2}} \left( \xi \right) \right) d \xi \left[\left[ \bold{V} \right]\right]_{i_{d}\pm\frac{1}{2}}
\end{equation}
Thus, 
\begin{equation}
\left(v^{(d)}_m\bold{F}_{m}\right)^{\star}_{i_{d}\pm\frac{1}{2}}=\frac{1}{2} \left( v^{(d)}_m\bold{F}_{m_i} + v^{(d)}_m\bold{F}_{m_{i_{d}\pm1}} \right) + O\left( \left| \left[\left[ \bold{V} \right]\right]_{i_{d}+\frac{1}{2}} \right|^2 \right)
\end{equation}
and hence for smooth functions, we have
\begin{equation}
\frac{1}{\Delta x_d} \left( \left(v^{(d)}_m\bold{F}_{m}\right)^{\star}_{i_{d}+\frac{1}{2}}-\left(v^{(d)}_m\bold{F}_{m}\right)^{\star}_{i_{d}-\frac{1}{2}}\right) =  \frac{1}{2 \Delta x_d} \left( \left(v^{(d)}_m\bold{F}_{m}\right)_{i_{d}+1}-\left(v^{(d)}_m\bold{F}_{m}\right)_{i_{d}-1}\right) + O \left( \left| \left[\left[ x_d \right]\right]_{i_{d}+\frac{1}{2}} \right|^2\right)
\end{equation}
Therefore, the entropy conserving scheme for vector-kinetic model given by \cref{avg+numdiff form for vec kin model} is second accurate in space. However, evaluation of a closed form interface flux function using \cref{avg+numdiff form for vec kin model} is algebraically tedious for a general hyperbolic system. \\
The closed form expression can be obtained along the same lines as macroscopic model in \cite{tadmor_2003}. Let $\left\{ \bold{l}^j_{i_{d}+\frac{1}{2}} \in \mathbb{R}^p \right\}_{j=1}^p$ and $\left\{ \bold{r}^j_{i_{d}+\frac{1}{2}} \in \mathbb{R}^p \right\}_{j=1}^p$ be two orthogonal sets of vectors such that $\left< \bold{l}^j_{i_{d}+\frac{1}{2}},  \bold{r}^k_{i_{d}+\frac{1}{2}} \right> = \delta_{jk}$. Let $\bold{V}^1_{i_{d}+\frac{1}{2}}=\bold{V}_i$ and
\begin{equation}
\bold{V}^{j+1}_{i_{d}+\frac{1}{2}} = \bold{V}^j_{i_{d}+\frac{1}{2}} + \left< \bold{l}^j_{i_{d}+\frac{1}{2}},  \left[\left[\bold{V}\right]\right]_{i_{d}+\frac{1}{2}} \right> \bold{r}^j_{i_{d}+\frac{1}{2}} \ ; j \in \{ 1,2,..,p\} 
\end{equation}
Then, we have a path connecting $\bold{V}_i$ and $\bold{V}_{i_{d}+1}$ since
\begin{equation}
\bold{V}^{p+1}_{i_{d}+\frac{1}{2}} = \bold{V}^1_{i_{d}+\frac{1}{2}} + \sum_{j=1}^p \left< \bold{l}^j_{i_{d}+\frac{1}{2}},  \left[\left[\bold{V}\right]\right]_{i_{d}+\frac{1}{2}} \right> \bold{r}^j_{i_{d}+\frac{1}{2}} = \bold{V}_i + \left[\left[\bold{V}\right]\right]_{i_{d}+\frac{1}{2}} = \bold{V}_{i_{d}+1}
\end{equation}
Now, it can be seen that the numerical flux given by, 
\begin{equation}
\label{PEC for vec kin model}
\left(v^{(d)}_m\bold{F}_{m}\right)^{\star}_{i_{d}+\frac{1}{2}}=\sum_{j=1}^p \frac{\chi^{(d)}_m\left( \bold{V}^{j+1}_{i_{d}+\frac{1}{2}} \right)- \chi^{(d)}_m\left( \bold{V}^{j}_{i_{d}+\frac{1}{2}}\right)}{\left< \bold{l}^j_{i_{d}+\frac{1}{2}},  \left[\left[\bold{V}\right]\right]_{i_{d}+\frac{1}{2}} \right>}  \bold{l}^j_{i_{d}+\frac{1}{2}}
\end{equation}
satisfies the entropy conserving condition in \cref{EC condn for vec kin model}. However, for the purpose of numerical simulations, we use robust entropy conserving fluxes (satisfying \cref{EC condn for vec kin model}) that are derived by defining averages of certain primitive variables and by balancing the coefficients corresponding to jumps in these primitive variables. These fluxes are described in \cref{Num res}. 
\begin{remark}
Higher order entropy conserving (HOEC) fluxes for vector-kinetic model can be constructed as linear combinations of second order entropy conserving fluxes derived in this paper (along the same lines as in \cite{doi:10.1137/S003614290240069X} for macroscopic model). Since linear combinations are used, as a consequence of \cref{thm EC for vec kin-mac model}, the moments of HOEC fluxes for vector-kinetic model will result in HOEC fluxes for macroscopic model.
\end{remark}
\begin{corollary}
\label{cor EC for vec-kin model}
If the assumptions stated in \cref{thm EC for vec kin-mac model} hold and entropy conserving flux of the form in \cref{avg+numdiff form for vec kin model} is used, then
\begin{equation}
\sum_{m=1}^M \bold{Q}^{(d)^{\star}}_{m_{i_{d}\pm\frac{1}{2}}}=\bold{Q}^{(d)^{\star}}_{i_{d}\pm\frac{1}{2}}
\end{equation}
\end{corollary}
\begin{proof}
By \cref{int flux mac model- vec kin model,avg+numdiff form for vec kin model}, we obtain
\begin{equation}
\bold{G}^{(d)^{\star}}_{i_{d}\pm\frac{1}{2}} =\sum_{m=1}^M \left(v^{(d)}_m\bold{F}_{m}\right)^{\star}_{i_{d}\pm\frac{1}{2}} = \frac{1}{2} \left( \bold{G}^{(d)}_i + \bold{G}^{(d)}_{i_{d}\pm1} \right) -\frac{1}{2} \sum_{m=1}^M \bold{Q}^{(d)^{\star}}_{m_{i_{d}\pm\frac{1}{2}}} \left[\left[ \bold{V} \right]\right]_{i_{d}\pm\frac{1}{2}}
\end{equation}
since $\sum_{m=1}^M v^{(d)}_m \bold{F}_{m_i}=\bold{G}^{(d)}_i, \ \forall i$ due to the action of moment constraint in \cref{Mom cons. vec kin model_2} on $\bold{F}_{m_i}$. Further,
\begin{equation}
\label{sum q for vec kin model}
\sum_{m=1}^M \bold{Q}^{(d)^{\star}}_{m_{i_{d}\pm\frac{1}{2}}} = \int_{0}^1 \left( 2 \xi -1 \right) \sum_{m=1}^M \partial_{\bold{V}} \left(v^{(d)}_m\bold{F}_{m}\right)\left( \bold{V}_{i_{d}+\frac{1}{2}} \left( \xi \right) \right) d \xi
\end{equation}
and
\begin{equation}
\sum_{m=1}^M \partial_{\bold{V}} \left(v^{(d)}_m\bold{F}_{m}\right)\left( \bold{V}_{i_{d}+\frac{1}{2}} \left( \xi \right) \right) = \sum_{m=1}^M v^{(d)}_m \partial_{\bold{V}} \left(  a_m \bold{U}  + b^j_m  \bold{G}^j \right) \left( \bold{V}_{i_{d}+\frac{1}{2}} \left( \xi \right) \right)= \partial_{\bold{V}} \bold{G}^{(d)}  \left( \bold{V}_{i_{d}+\frac{1}{2}} \left( \xi \right) \right)
\end{equation}
due to the action of moment constraint in \cref{Mom cons. vec kin model_2} on $\partial_{\bold{V}}\bold{F}_m$. Thus, comparing \cref{sum q for vec kin model,Q for mac model}, we obtain $\sum_{m=1}^M \bold{Q}^{(d)^{\star}}_{m_{i_{d}\pm\frac{1}{2}}}=\bold{Q}^{(d)^{\star}}_{i_{d}\pm\frac{1}{2}}$. 
\end{proof}

\section{Entropy stable scheme for vector-kinetic model}
\label{Sec:ES}
Consider the three-point semi-discrete conservative scheme on structured grid,
\begin{equation}
\label{ES scheme for vec kin model}
\frac{d}{d t} \bold{F}_{m_i} + \frac{1}{\Delta x_d} \left( \left(v^{(d)}_m\bold{F}_{m}\right)_{i_{d}+\frac{1}{2}}-\left(v^{(d)}_m\bold{F}_{m}\right)_{i_{d}-\frac{1}{2}}\right)=\bold{0}
\end{equation}
The interface numerical flux $\left(v^{(d)}_m\bold{F}_{m}\right)_{i_{d}\pm\frac{1}{2}}$ is given by,
\begin{equation} 
\label{ES int flux for vec kin model}
\left(v^{(d)}_m\bold{F}_{m}\right)_{i_{d}\pm\frac{1}{2}}= \left(v^{(d)}_m\bold{F}_{m}\right)^{\star}_{i_{d}\pm\frac{1}{2}}-\frac{1}{2} \bold{D}^{(d)}_{m_{i_{d}\pm\frac{1}{2}}} \left[\left[ \partial_{\bold{F}_m} H^{\eta}_m \right]\right]_{i_{d}\pm\frac{1}{2}}
\end{equation}
Here, $\bold{D}^{(d)}_{m_{i_{d}\pm\frac{1}{2}}}=\bold{Q}^{(d)}_{m_{i_{d}\pm\frac{1}{2}}}-\bold{Q}^{(d)^{\star}}_{m_{i_{d}\pm\frac{1}{2}}}$. $\bold{Q}^{(d)}_{m_{i_{d}\pm\frac{1}{2}}}$ and $\bold{Q}^{(d)^{\star}}_{m_{i_{d}\pm\frac{1}{2}}}$ are the numerical viscosity coefficient matrices corresponding to entropy stable and entropy conserving schemes respectively. $\bold{Q}^{(d)^{\star}}_{m_{i_{d}\pm\frac{1}{2}}}$ is given by \cref{num visc of ec_kin}. \\
Then, the inner product of \cref{ES scheme for vec kin model} with $\left(\partial_{\bold{F}_m} H^{\eta}_m\right)_i$ gives the entropy in-equality,
\begin{multline}
\label{ES eqn for vec kin model}
\frac{d}{d t} H^{\eta}_{m_i}  +\frac{1}{\Delta x_d} \left( \left(v^{(d)}_mH^{\eta}_{m}\right)_{i_{d}+\frac{1}{2}}-\left(v^{(d)}_mH^{\eta}_{m}\right)_{i_{d}-\frac{1}{2}}\right)\\ = -\frac{1}{4\Delta x_d}  \left( \left[\left[ \partial_{\bold{F}_m} H^{\eta}_m \right]\right]_{i_{d}+\frac{1}{2}} \cdot \bold{D}^{(d)}_{m_{i_{d}+\frac{1}{2}}} \left[\left[ \partial_{\bold{F}_m} H^{\eta}_m \right]\right]_{i_{d}+\frac{1}{2}} + \left[\left[ \partial_{\bold{F}_m} H^{\eta}_m \right]\right]_{i_{d}-\frac{1}{2}} \cdot \bold{D}^{(d)}_{m_{i_{d}-\frac{1}{2}}} \left[\left[ \partial_{\bold{F}_m} H^{\eta}_m \right]\right]_{i_{d}-\frac{1}{2}} \right) \leq 0
\end{multline}
iff $\bold{D}^{(d)}_{m_{i_{d}\pm\frac{1}{2}}}$ is positive-definite. The interface numerical entropy flux $\left(v^{(d)}_mH^{\eta}_{m}\right)_{i_{d}+\frac{1}{2}}$ consistent with \cref{efp for vec kin model} becomes, 
\begin{equation}
\label{ES Ent int flux for vec kin model}
 \left(v^{(d)}_mH^{\eta}_{m}\right)_{i_{d}+\frac{1}{2}}= \left(v^{(d)}_mH^{\eta}_{m}\right)^{\star}_{i_{d}+\frac{1}{2}} - \frac{1}{4} \left( \left( \partial_{\bold{F}_m} H^{\eta}_m\right)_i+\left(\partial_{\bold{F}_m} H^{\eta}_m \right)_{i_{d}+1} \right) \cdot \bold{D}^{(d)}_{m_{i_{d}+\frac{1}{2}}} \left[\left[ \partial_{\bold{F}_m} H^{\eta}_m \right]\right]_{i_{d}+\frac{1}{2}} 
\end{equation}
Further, let us define the interface numerical fluxes for macroscopic model as the moment of interface numerical fluxes for vector-kinetic model as,
\begin{equation}
\label{ES int flux mac model- vec kin model}
\bold{G}^{(d)}_{i_{d}\pm\frac{1}{2}} = \sum_{m=1}^M \left(v^{(d)}_m\bold{F}_{m}\right)_{i_{d}\pm\frac{1}{2}}
\end{equation}
\begin{theorem}
\label{thm ES for vec kin-mac model}
If the three-point semi-discrete conservative scheme $\left(\cref{ES scheme for vec kin model}\right)$ for vector-kinetic model with 
\begin{itemize}
\item{$\bold{F}_{m_i} =a_m\bold{U}_i+b^{(d)}_m\bold{G}^{(d)}_i, \ \forall i$}
\item{interface numerical fluxes $\left(v^{(d)}_m\bold{F}_{m}\right)_{i_{d}\pm\frac{1}{2}}$ satisfying \cref{ES int flux for vec kin model} and} 
\item{constants $a_m$, $b^{(d)}_m$ satisfying the moment constraints in \cref{Mom cons. vec kin model_1,Mom cons. vec kin model_2} while rendering the positivity of eigenvalues of $\partial_{\bold{U}} \bold{F}_m $}
\end{itemize}
is used, and if the convex entropy function corresponding to it is $H^{\eta}_{m_i} =a_m\eta_i+b^{(d)}_m\omega^{(d)}_i, \ \forall i$, then 
\begin{enumerate}
\item{$\sum_{m=1}^M$ \cref{ES scheme for vec kin model} becomes 
\begin{equation}
\label{thm ES 1st item}
\frac{d}{d t} \bold{U}_i + \frac{1}{\Delta x_d} \left( \bold{G}^{(d)}_{i_{d}+\frac{1}{2}}-\bold{G}^{(d)}_{i_{d}-\frac{1}{2}}\right)=\bold{0}
\end{equation} 
with $\bold{G}^{(d)}_{i_{d}\pm\frac{1}{2}}$ given by \cref{ES int flux mac model- vec kin model},}
\item{the interface numerical flux $\bold{G}^{(d)}_{i_{d}\pm\frac{1}{2}}$ given by \cref{ES int flux mac model- vec kin model} is equal to \cref{ES flux for mac model}, and}
\item{the three-point entropy in-equality for macroscopic model $\left(\cref{ES eqn for mac model}\right)$ holds true with interface numerical entropy flux $\omega^{(d)}_{i_{d}\pm\frac{1}{2}}$ given by \cref{ES Ent flux for mac model}.}
\end{enumerate}
\end{theorem}
\begin{proof}
Due to moment constraint in \cref{Mom cons. vec kin model_1}, $\sum_{m=1}^M \bold{F}_{m_i}=\bold{U}_i$. Therefore, $\sum_{m=1}^M$ \cref{ES scheme for vec kin model} becomes \cref{thm ES 1st item} with $\bold{G}^{(d)}_{i_{d}\pm\frac{1}{2}}$ given by \cref{ES int flux mac model- vec kin model}, thus proving 1. \\
Since $\left(v^{(d)}_m\bold{F}_{m}\right)_{i_{d}\pm\frac{1}{2}}$ follows \cref{ES int flux for vec kin model} and $\left[\left[\partial_{\bold{F}_m} H^{\eta}_m \right]\right]_{i_{d}\pm\frac{1}{2}}=\left[\left[\bold{V}\right]\right]_{i_{d}\pm\frac{1}{2}}=  \left[\left[\partial_{\bold{U}}  \eta\right]\right]_{i_{d}\pm\frac{1}{2}}$ is not a function of $m$ (by \cref{equality of ent var}), \cref{ES int flux mac model- vec kin model} becomes,
\begin{equation}
\bold{G}^{(d)}_{i_{d}\pm\frac{1}{2}} = \sum_{m=1}^M \left(v^{(d)}_m\bold{F}_{m}\right)_{i_{d}\pm\frac{1}{2}} = \sum_{m=1}^M \left(v^{(d)}_m\bold{F}_{m}\right)^{\star}_{i_{d}\pm\frac{1}{2}}-\frac{1}{2} \sum_{m=1}^M \bold{D}^{(d)}_{m_{i_{d}\pm\frac{1}{2}}} \left[\left[ \bold{V} \right]\right]_{i_{d}\pm\frac{1}{2}}
\end{equation}
By \cref{thm EC for vec kin-mac model}, $\sum_{m=1}^M \left(v^{(d)}_m\bold{F}_{m}\right)^{\star}_{i_{d}\pm\frac{1}{2}}$ satisfies entropy conserving condition in \cref{EC condn for mac model} and hence it is equal to $\bold{G}^{(d)^{\star}}_{i_{d}\pm\frac{1}{2}}$. We also have $\sum_{m=1}^M \bold{Q}^{(d)^{\star}}_{m_{i_{d}\pm\frac{1}{2}}}=\bold{Q}^{(d)^{\star}}_{i_{d}\pm\frac{1}{2}}$ by \cref{cor EC for vec-kin model}. Further, $\sum_{m=1}^M \bold{D}^{(d)}_{m_{i_{d}\pm\frac{1}{2}}}$ is positive-definite as $\bold{D}^{(d)}_{m_{i_{d}\pm\frac{1}{2}}}$ is positive-definite $\forall m$. Therefore, $\bold{D}^{(d)}_{i_{d}\pm\frac{1}{2}} = \sum_{m=1}^M \bold{D}^{(d)}_{m_{i_{d}\pm\frac{1}{2}}} = \sum_{m=1}^M \bold{Q}^{(d)}_{m_{i_{d}\pm\frac{1}{2}}}- \bold{Q}^{(d)^{\star}}_{i_{d}\pm\frac{1}{2}}$ is positive-definite, and hence
\begin{equation}
\bold{G}^{(d)}_{i_{d}\pm\frac{1}{2}}= \bold{G}^{(d)^{\star}}_{i_{d}\pm\frac{1}{2}} -\frac{1}{2} \bold{D}^{(d)}_{i_{d}\pm\frac{1}{2}} \left[\left[ \bold{V} \right]\right]_{i_{d}\pm\frac{1}{2}}
\end{equation}
This proves 2. \\
Corresponding to the assumptions stated in \cref{thm ES for vec kin-mac model}, the three-point entropy in-equality in \cref{ES eqn for vec kin model} holds true. Since $\sum_{m=1}^M (H^{\eta}_m)_i=\eta_i, \ \forall i$ (due to the action of moment constraint in \cref{Mom cons. vec kin model_1} on $(H^{\eta}_m)_i$), $\left[\left[\partial_{\bold{F}_m} H^{\eta}_m \right]\right]_{i_{d}\pm\frac{1}{2}}=\left[\left[\bold{V}\right]\right]_{i_{d}\pm\frac{1}{2}}=  \left[\left[\partial_{\bold{U}}  \eta\right]\right]_{i_{d}\pm\frac{1}{2}}$ is not a function of $m$ (by \cref{equality of ent var}) and $\sum_{m=1}^M \bold{D}^{(d)}_{m_{i_{d}+\frac{1}{2}}} = \bold{D}^{(d)}_{i_{d}+\frac{1}{2}}$, moment of \cref{ES eqn for vec kin model} gives,
\begin{multline}
\frac{d}{d t} \eta_i + \frac{1}{\Delta x_d} \left( \sum_{m=1}^M \left(v^{(d)}_mH^{\eta}_{m}\right)_{i_{d}+\frac{1}{2}}-\sum_{m=1}^M \left(v^{(d)}_mH^{\eta}_{m}\right)_{i_{d}-\frac{1}{2}}\right) = \\ -\frac{1}{4\Delta x_d}  \left( \left[\left[ \bold{V} \right]\right]_{i_{d}+\frac{1}{2}} \cdot  \bold{D}^{(d)}_{i_{d}+\frac{1}{2}} \left[\left[ \bold{V} \right]\right]_{i_{d}+\frac{1}{2}} + \left[\left[ \bold{V} \right]\right]_{i_{d}-\frac{1}{2}} \cdot \bold{D}^{(d)}_{i_{d}-\frac{1}{2}} \left[\left[ \bold{V} \right]\right]_{i_{d}-\frac{1}{2}} \right)
\end{multline}
Since $\left[\left[\partial_{\bold{F}_m} H^{\eta}_m \right]\right]_{i_{d}\pm\frac{1}{2}}=\left[\left[\bold{V}\right]\right]_{i_{d}\pm\frac{1}{2}}=  \left[\left[\partial_{\bold{U}}  \eta\right]\right]_{i_{d}\pm\frac{1}{2}}$ and $\left(\partial_{\bold{F}_m} H^{\eta}_m\right)_i=\bold{V}_i=\left(\partial_{\bold{U}}  \eta\right)_i$ are not functions of $m$ (by \cref{equality of ent var}), moment of \cref{ES Ent int flux for vec kin model} yields,
\begin{equation}
\sum_{m=1}^M \left(v^{(d)}_mH^{\eta}_{m}\right)_{i_{d}+\frac{1}{2}} =  \sum_{m=1}^M \left(v^{(d)}_mH^{\eta}_{m}\right)^{\star}_{i_{d}+\frac{1}{2}} - \frac{1}{4} \left( \bold{V}_i+\bold{V}_{i_{d}+1} \right) . \sum_{m=1}^M \bold{D}^{(d)}_{m_{i_{d}+\frac{1}{2}}} \left[\left[ \bold{V} \right]\right]_{i_{d}+\frac{1}{2}} 
\end{equation}
Since  $\sum_{m=1}^M \left(v^{(d)}_mH^{\eta}_{m}\right)^{\star}_{i_{d}+\frac{1}{2}}=\omega^{(d)^{\star}}_{i_{d}+\frac{1}{2}}$ (by \cref{thm EC for vec kin-mac model}) and $\sum_{m=1}^M \bold{D}^{(d)}_{m_{i_{d}+\frac{1}{2}}} = \bold{D}^{(d)}_{i_{d}+\frac{1}{2}}$, comparison of the above equation with \cref{ES Ent flux for mac model} yields $\sum_{m=1}^M \left(v^{(d)}_mH^{\eta}_{m}\right)_{i_{d}+\frac{1}{2}} =\omega^{(d)}_{i_{d}+\frac{1}{2}}$. This proves 3. 
\end{proof}
Thus, an entropy stable scheme for vector-kinetic model results in an entropy stable scheme for macroscopic model, thanks to \cref{equality of ent var} (resulting from \cref{Exis ent var for vec kin model}) that rendered the linearity of moments in the above proof. 
\begin{remark}
In the above proof, the three-point entropy in-equality for macroscopic model $\left(\cref{ES eqn for mac model}\right)$ with interface numerical entropy flux $\omega^{(d)}_{i_{d}\pm\frac{1}{2}}$ given by \cref{ES Ent flux for mac model} is obtained as moment of three-point entropy in-equality for vector-kinetic model. Unlike this, we can also obtain \cref{ES eqn for mac model} directly at the macroscopic level as a consequence of $\bold{G}^{(d)}_{i_{d}\pm\frac{1}{2}}= \sum_{m=1}^M \left(v^{(d)}_m\bold{F}_{m}\right)_{i_{d}\pm\frac{1}{2}}$ satisfying the entropy stability condition for macroscopic model $\left(\cref{ES flux for mac model} \text{ with positive-definite } \bold{D}^{(d)}_{i_{d}\pm\frac{1}{2}}\right)$.  
\end{remark}

\subsection{High resolution scheme}
Since the interface numerical flux $\left(v^{(d)}_m\bold{F}_{m}\right)_{i_{d}+\frac{1}{2}}$ contains a term with $\left[\left[ \bold{V} \right]\right]_{i_{d}+\frac{1}{2}}$ which is $O\left( \Delta x_d \right)$, the entropy stable scheme in \cref{ES scheme for vec kin model} is only first order accurate in space. In order to attain higher order accuracy in space, the interface numerical flux in \cref{ES int flux for vec kin model} is modified as,
\begin{equation} 
\label{HO ES int flux for vec kin model}
\left(v^{(d)}_m\bold{F}_{m}\right)_{i_{d}\pm\frac{1}{2}}= \left(v^{(d)}_m\bold{F}_{m}\right)^{\star}_{i_{d}\pm\frac{1}{2}}-\frac{1}{2} \bold{D}^{(d)}_{m_{i_{d}\pm\frac{1}{2}}} \left<\left< \bold{V} \right>\right>_{i_{d}\pm\frac{1}{2}}
\end{equation}
where $\left<\left< \bold{V} \right>\right>_{i_{d}+\frac{1}{2}}=\bold{V}^-_{i_{d}+1}-\bold{V}^+_i$. Further, $\bold{V}^-_{i_{d}+1}=\bold{V}_{i_{d}+1}\left(x_{d_{i_{d}+\frac{1}{2}}}\right)$ and $\bold{V}^+_i=\bold{V}_{i}\left(x_{d_{i_{d}+\frac{1}{2}}}\right)$ are higher order reconstructions of $\bold{V}$ at interface $i_{d}+\frac{1}{2}$. We utilise second order reconstructions in obtaining the numerical results, and the details are provided therein \cref{Sec:NR}. The moment of \cref{HO ES int flux for vec kin model} becomes, 
\begin{equation}
\sum_{m=1}^M \left(v^{(d)}_m\bold{F}_{m}\right)_{i_{d}\pm\frac{1}{2}}= \bold{G}^{(d)^{\star}}_{i_{d}\pm\frac{1}{2}}-\frac{1}{2} \bold{D}^{(d)}_{i_{d}\pm\frac{1}{2}} \left<\left< \bold{V} \right>\right>_{i_{d}\pm\frac{1}{2}}
\end{equation}
It can be easily seen that this is a higher order entropy stable flux for macroscopic model, and it is a consequence of linearity due to \cref{equality of ent var} (resulting from \cref{Exis ent var for vec kin model}).  

\section{Time discretisation}
\label{Sec:TD}
Let $\mathcal{F}_{m_i}$ be $-\frac{1}{\Delta x_d} \left( \left(v^{(d)}_m\bold{F}_{m}\right)_{i_{d}+\frac{1}{2}}-\left(v^{(d)}_m\bold{F}_{m}\right)_{i_{d}-\frac{1}{2}}\right)$ where $\left(v^{(d)}_m\bold{F}_{m}\right)_{i_{d}\pm\frac{1}{2}}$ is entropy conserving $\biggl(\biggr.$ $\left(v^{(d)}_m\bold{F}_{m}\right)^{\star}_{i_{d}\pm\frac{1}{2}}$ satisfying \cref{EC condn for vec kin model}$\biggl.\biggr)$ or entropy stable $\biggl(\biggr.$ $\left(v^{(d)}_m\bold{F}_{m}\right)_{i_{d}\pm\frac{1}{2}}$ satisfying \cref{ES int flux for vec kin model}$\biggl. \biggr)$. Then, the semi-discrete entropy conserving/stable schemes in \cref{EC scheme for vec kin model,ES scheme for vec kin model} can be represented as,
\begin{equation}
\frac{d}{d t} \bold{F}_{m_i} = \mathcal{F}_{m_i}
\end{equation}
Since we utilise second order scheme for entropy conserving/stable spatial discretisations, a third order scheme is required for the temporal derivative so that the entropy production/dissipation due to temporal derivative will not affect the entropy conservation/stability achieved spatially. Hence, the temporal derivative in above equation is discretised using 3-stage third order strong stability preserving Runge-Kutta method ($\text{SSPRK}(3,3)$) \cite{SHU1988439}. After each stage of the RK method, $\bold{U}_i$ is evaluated using $\bold{U}_i=\sum_{m=1}^M \bold{F}_{m_i}$, and this is utilised in the evaluation of fluxes required for the next stage. 

\section{Numerical results}
\label{Sec:NR}
\label{Num res}
In this section, the entropy conserving (EC)/stable (ES) schemes are tested against various physical problems governed by scalar equations and the system of shallow water equations. For each problem, the basic ingredients such as problem description, choice of macroscopic entropy-entropy flux pair, fluxes satisfying entropy conserving/stability conditions in \cref{EC condn for vec kin model,ES int flux for vec kin model}, second order reconstructions of entropy stable fluxes and CFL criteria are provided. We use the following error quantifications to study the errors in macroscopic and vector-kinetic entropies at time $t$. 
\begin{gather}
\text{Signed error} = \frac{ \sum_i \left((.)^{t}_i - (.)^{t-\Delta t}_i \right) }{N} \\
\text{Absolute error} =  \frac{ \sum_i \left| (.)^{t}_i - (.)^{t-\Delta t}_i  \right| }{N}
\end{gather}
Here, $N$ is the total number of cells or grid points in the computational domain. It can be seen that the signed error allows for cancellations of positive and negative errors present at different spatial locations. An equivalent of this with reference as $t=0$ instead of $t-\Delta t$ is commonly used in literature in the context of global entropy preservation \cite{doi:10.1137/19M1263480}. However, in order to understand the actual entropy preservation property of a spatially entropy preserving scheme, one needs to use the absolute error that does not allow spatial cancellations. Further, we use the signed error to identify whether the scheme is globally entropy dissipating or not. A positive signed error indicates global entropy production while negative signed error indicates global entropy dissipation. We present the numerical solutions, global entropy vs. time, and error vs. time plots for each problem.

\subsection{Scalar equations}
We consider scalar equations of the form,
\begin{equation}
\partial_t U + \partial_{x_d} G^{(d)}(U) = 0
\end{equation}
with initial condition $U(x_1,..,x_d,..,x_D,0)=U_0(x_1,..,x_d,..,x_D)$. We choose suitable convex entropy-entropy flux pair specific to $G^{(d)}(U)$. The constants $a_m, b^{(d)}_m$ in \cref{F def vec kin model,H def vec kin model} are chosen as described in \ref{app constants a,b}. The time step is chosen as
\begin{equation}
\Delta t \leq C \frac{\Delta x }{\lambda} \ ; \  \Delta x = min\left( \Delta x_d \right)
\end{equation}
Here, $C$ is the CFL number. The choice of $\lambda$ is described in \ref{app constants a,b}. The flux
\begin{equation}
\left(v^{(d)}_mF_{m}\right)^{\star}_{i_{d}+\frac{1}{2}} = \frac{\chi^{(d)}_{m_{i_{d}+1}}-\chi^{(d)}_{m_{i}}}{V_{i_{d}+1}-V_i}
\end{equation}
satisfies the entropy conserving condition in \cref{EC condn for vec kin model}. This is used when $V_{i_{d}+1}\neq V_i$. When $V_{i_{d}+1} =V_i$, we do not update the flux, as any value of flux satisfies the entropy conserving condition $\left(\cref{EC condn for vec kin model}\right)$. Here, the entropy variable is $V_i= \left( \partial_{U} \eta \right)_i$ and the vector-kinetic entropy flux potential is given by $\chi^{(d)}_{m_i}=V_i . \left( v^{(d)}_m F_m\right)_i - \left( v^{(d)}_m H^{\eta}_m\right)_i$. \\
For entropy stable scheme, we use $\bold{D}^{(d)}_{m_{i_{d}+\frac{1}{2}}} \left<\left< \bold{V}\right>\right>_{i_{d}+\frac{1}{2}} = \frac{1}{M}\bold{R}^{(d)}_{i_{d}+\frac{1}{2}} \bold{\Lambda}^{(d)}_{i_{d}+\frac{1}{2}} \left<\left< \widetilde{\bold{W}}\right>\right>_{i_{d}+\frac{1}{2}}$. For scalar equations, $\bold{R}^{(d)}_{i_{d}+\frac{1}{2}}=1$ and $\bold{\Lambda}^{(d)}_{i_{d}+\frac{1}{2}}$ is the absolute wave speed obtained using the average (arithmetic) value of $U$ at cells $i$ and $i_{d}+1$. We use the second order reconstruction of $\left<\left< \widetilde{\bold{W}}\right>\right>_{i_{d}+\frac{1}{2}}$ as explained in \cref{SubSec: SW}. 

\subsubsection{Linear advection}
\begin{figure}[t]
\centering
\begin{subfigure}[b]{0.23\textwidth}
\centering
\includegraphics[width=\textwidth]{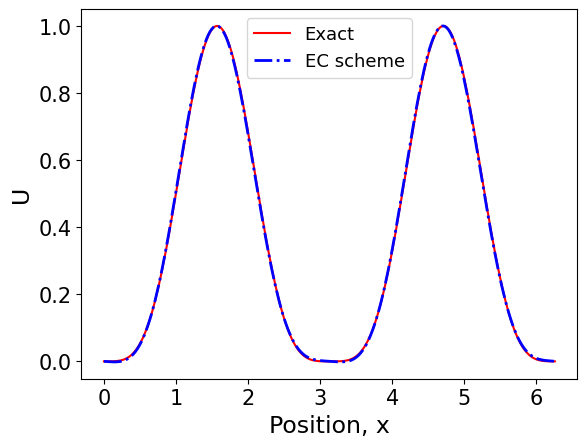}
\caption{Solutions}
\label{sc la sol}
\end{subfigure}
\hfill
\begin{subfigure}[b]{0.23\textwidth}
\centering
\includegraphics[width=\textwidth]{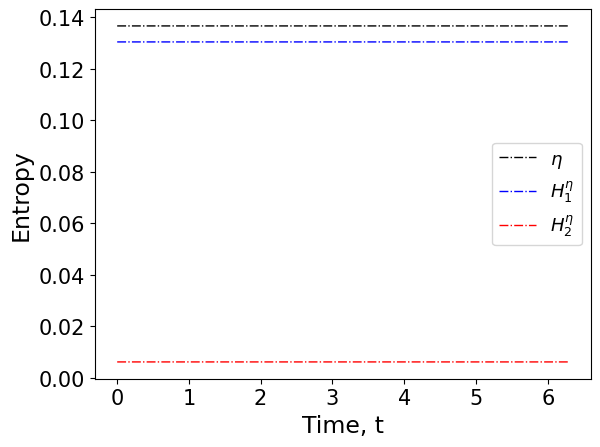}
\caption{Entropy functions}
\label{sc la ent}
\end{subfigure}
\hfill 
\begin{subfigure}[b]{0.23\textwidth}
\centering
\includegraphics[width=\textwidth]{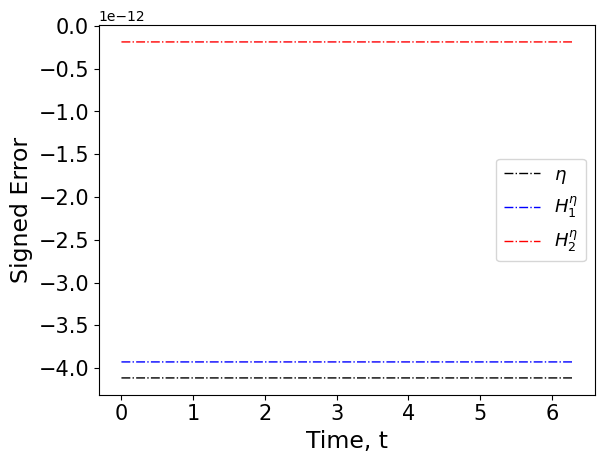}
\caption{Signed errors}
\label{sc la Serr}
\end{subfigure}
\hfill
\begin{subfigure}[b]{0.24\textwidth}
\centering
\includegraphics[width=\textwidth]{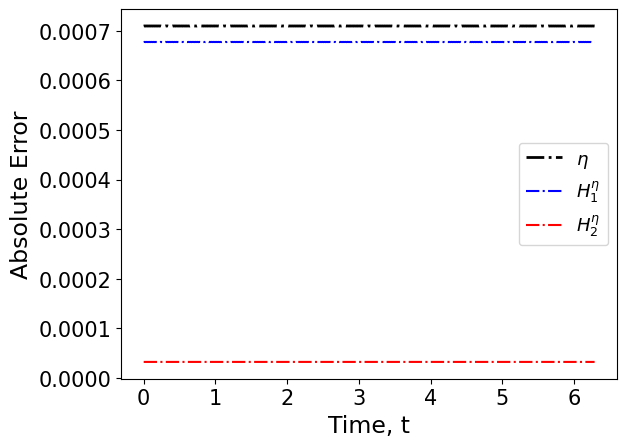}
\caption{Absolute errors}
\label{sc la USerr}
\end{subfigure}
\caption{\centering Linear advection at $T=2\pi$ using EC scheme with $C=0.1$ and $Nx=256$}
\label{sc la}
\end{figure}

For the one-dimensional linear advection problem with $G^{(1)}(U)=U$, we choose $\eta(U)=\frac{1}{2}U^2$, and correspondingly $\omega^{(1)}(U)=\frac{1}{2}U^2$ satisfies the compatibility condition in \cref{Ent codn for mac model}.  The initial condition is $U_0\left(x_1\right)=\left(sin (x_1)\right)^4$. The domain of the problem is $[0,2\pi)$, and it is discretised using 256 uniform cells. Periodic boundary conditions are used here. Numerical solutions are obtained at $T=2\pi$.\\    
It can be seen from \cref{sc la sol} that the numerical solution matches well with the exact solution. \Cref{sc la ent} shows the global entropies over time. It can be seen that the entropies remain nearly constant. The signed and absolute errors in entropies are shown in \cref{sc la Serr,sc la USerr} respectively. Since we use second order accurate entropy conserving scheme for vector-kinetic model and $\Delta x$ is of $O(10^{-2})$, we expect an absolute error of $O(10^{-4})$ in the vector-kinetic entropies. This is observed in \cref{sc la USerr}. The negative signed errors in \cref{sc la Serr} indicate that the $O(\Delta x^2)$ error is globally dissipative in nature. Due to the symmetric nature of the periodic profile, there may be cancellations in errors spatially and we observe a very low signed error of $O(10^{-12})$. \\
In order to study the convergence of the problem, we use very low CFL of $C=0.1$. Second order accuracy of the scheme is evident from the results presented in \cref{tab:EOC_la}. The exact solution is used as reference for the convergence study. 
\begin{table}[tbhp]
\begin{center}
\begin{tabular}{|m{2cm}|m{2cm}|m{2cm}|m{2cm}|}
\hline
\centering Number of cells, Nx & \centering $\Delta x_1$ & \centering $L_2$ norm & $O(L_2)$ \\
\hline
\centering 32 & \centering 0.196349541 & \centering 0.035757668 &  - \\
\centering 64 & \centering 0.09817477 & \centering 0.00781911 &  2.19 \\
\centering 128 & \centering 0.049087385 & \centering 0.00140703 &  2.47 \\
\centering 256 & \centering 0.024543693 & \centering 0.000249239 &  2.50 \\
\hline
\end{tabular}
\caption{\centering EOC for linear advection at $T=2\pi$ using EC scheme with $C=0.1$} 
\label{tab:EOC_la}
\end{center}
\end{table}

\subsubsection{Linear rotation}
\begin{figure}[t]
\centering
\begin{subfigure}[b]{0.3\textwidth}
\centering
\includegraphics[width=\textwidth]{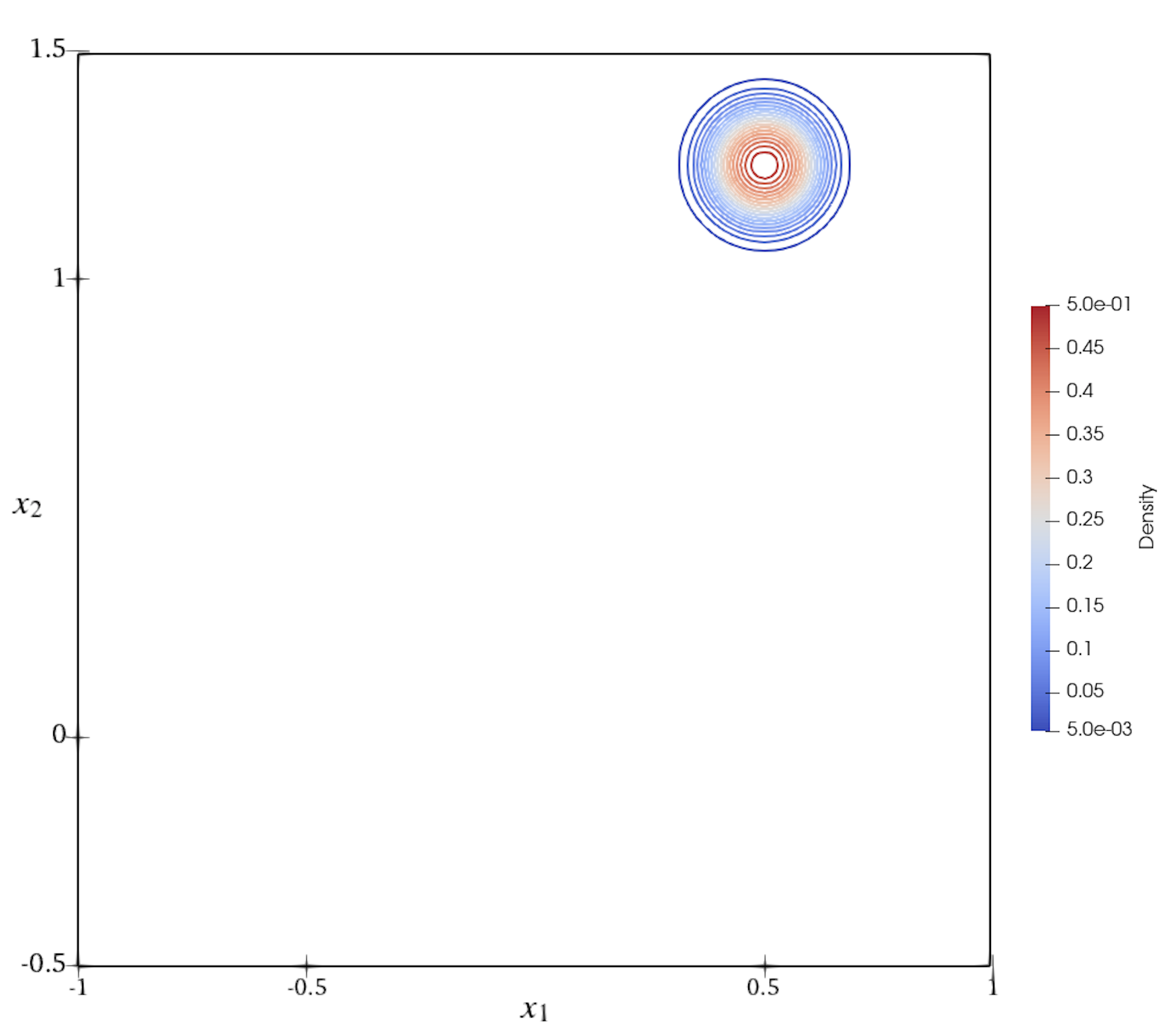}
\caption{Initial condn. $T=0$}
\label{sc lr ic}
\end{subfigure}
\hfill
\begin{subfigure}[b]{0.3\textwidth}
\centering
\includegraphics[width=\textwidth]{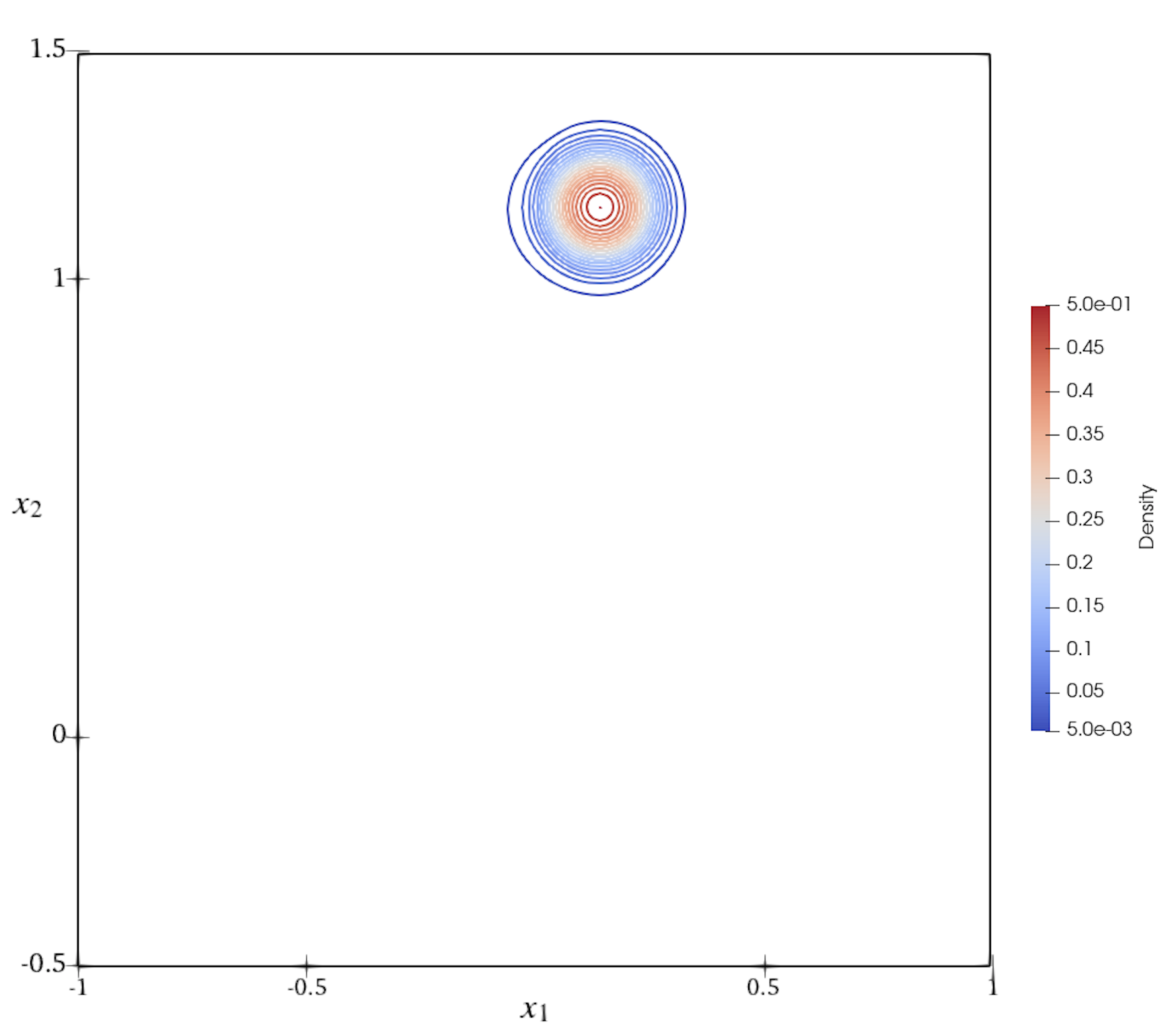}
\caption{Numerical soln. $T=0.5$}
\label{sc lr sol}
\end{subfigure}
\vfill
\begin{subfigure}[b]{0.28\textwidth}
\centering
\includegraphics[width=\textwidth]{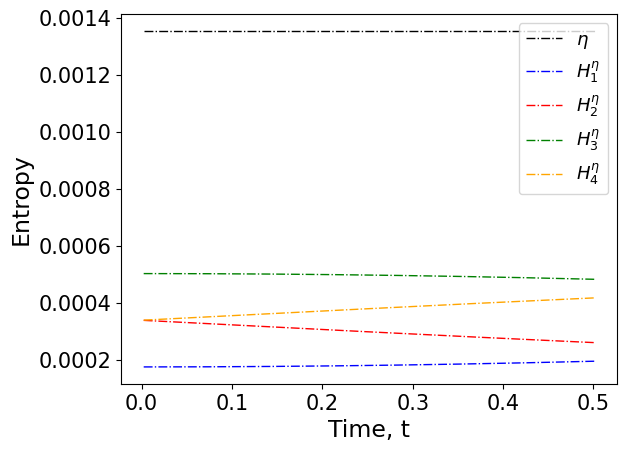}
\caption{Entropy functions}
\label{sc lr ent}
\end{subfigure}
\hfill 
\begin{subfigure}[b]{0.26\textwidth}
\centering
\includegraphics[width=\textwidth]{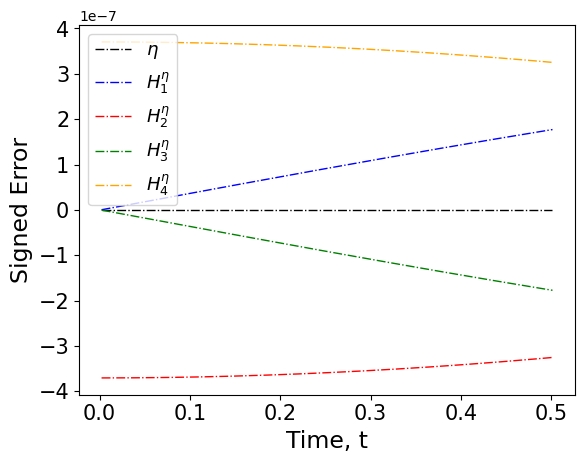}
\caption{Signed errors}
\label{sc lr Serr}
\end{subfigure}
\hfill 
\begin{subfigure}[b]{0.26\textwidth}
\centering
\includegraphics[width=\textwidth]{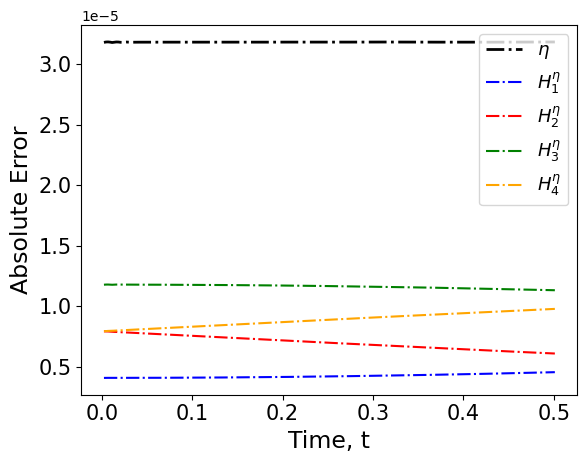}
\caption{Absolute errors}
\label{sc lr USerr}
\end{subfigure}
\caption{\centering Linear rotation at $T=0.5$ using EC scheme with $C=0.9$ and $Nx, Ny=256$}
\label{sc lr Err}
\end{figure}

For the two dimensional linear rotation problem, $G^{(1)}(U)=-\left(x_2 - \frac{1}{2}\right)U$ and $G^{(2)}(U)=\left(x_1 - \frac{1}{2}\right)U$. The entropy function is chosen as $\eta(U)=U^2$, and correspondingly the entropy flux functions become $\omega^{(1)}(U)=-\left(x_2 - \frac{1}{2}\right)U^2$ and $\omega^{(2)}(U)=\left(x_1 - \frac{1}{2}\right)U^2$. The initial condition is shown in \cref{sc lr ic}. The domain of the problem is $[-1,1) \times [-0.5,1.5)$, and it is discretised using $256 \times 256$ uniform cells. The value of $U$ at the boundary is kept fixed throughout the computation, and a CFL of $C=0.9$ is used. \\
The numerical solution at $T=0.5$ is shown in \cref{sc lr sol}. Since $\Delta x$ is of $O(10^{-2})$, one would expect an error of $O(10^{-4})$ in the absolute errors due to the usage of second order accurate entropy conserving scheme. We observe better error of $O(10^{-5})$ in \cref{sc lr USerr}. Further, it is interesting to observe the symmetries in errors of $H^{\eta}_2, H^{\eta}_4$ and $H^{\eta}_1, H^{\eta}_3$ in \cref{sc lr Serr}. However, these symmetries may not be located on the same spatial point. If they were, then the absolute error of macroscopic entropy $\eta$ would be much smaller than $O(10^{-7})$ (due to cancellations) since it is the sum of vector-kinetic entropies.

\subsubsection{Non-linear inviscid Burgers' test}
\begin{figure}[t]
\centering
\begin{subfigure}[b]{0.23\textwidth}
\centering
\includegraphics[width=\textwidth]{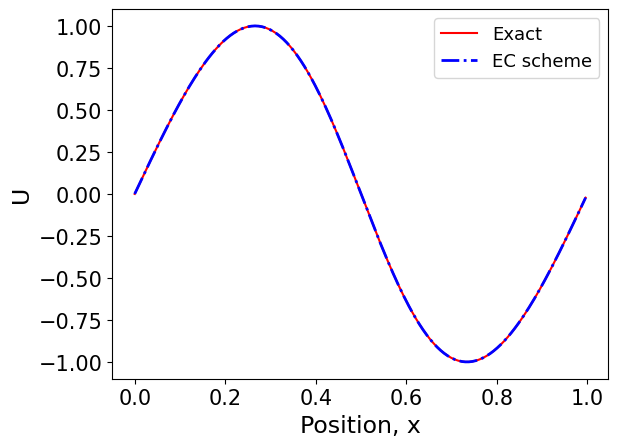}
\caption{Solutions}
\label{sc nlib sol}
\end{subfigure}
\hfill
\begin{subfigure}[b]{0.23\textwidth}
\centering
\includegraphics[width=\textwidth]{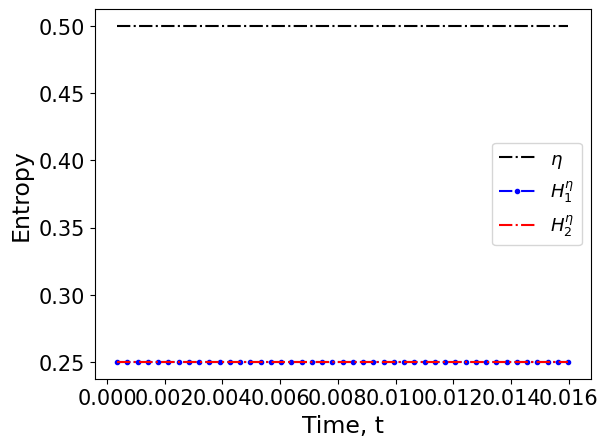}
\caption{Entropy functions}
\label{sc nlib ent}
\end{subfigure}
\hfill 
\begin{subfigure}[b]{0.23\textwidth}
\centering
\includegraphics[width=\textwidth]{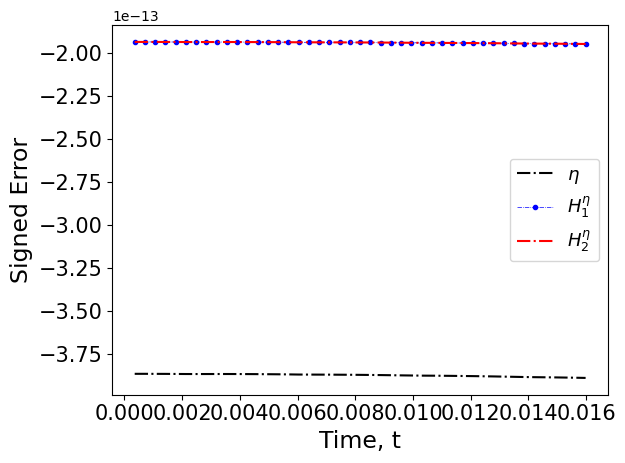}
\caption{Signed errors}
\label{sc nlib Serr}
\end{subfigure}
\hfill
\begin{subfigure}[b]{0.23\textwidth}
\centering
\includegraphics[width=\textwidth]{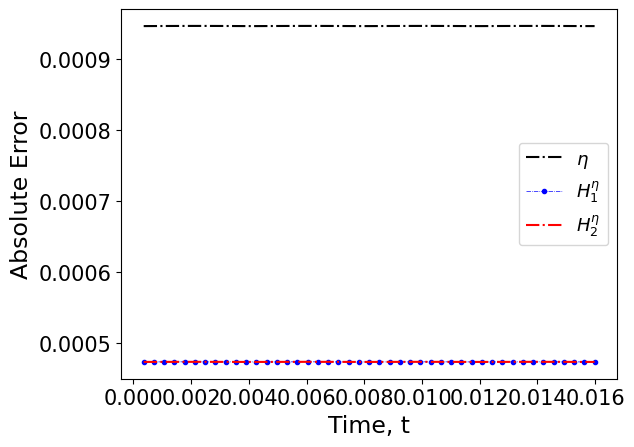}
\caption{Absolute errors}
\label{sc nlib USerr}
\end{subfigure}
\caption{\centering Non-linear inviscid Burgers' test at $T=\frac{0.1}{2\pi}$ using EC scheme with $C=0.1$ and $Nx=256$}
\label{sc nlib}
\end{figure}

\begin{figure}
\centering
\begin{subfigure}[b]{0.23\textwidth}
\centering
\includegraphics[width=\textwidth]{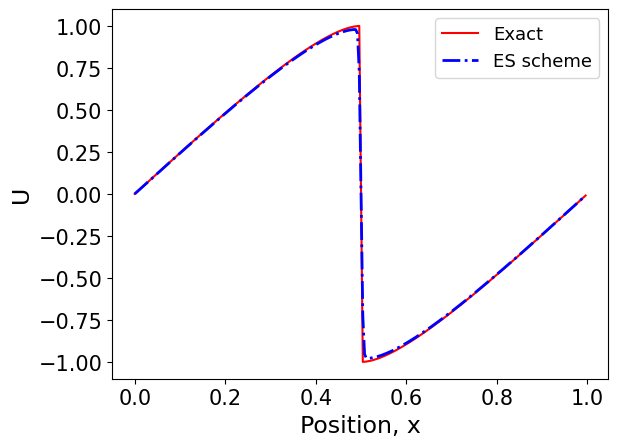}
\caption{Solutions}
\label{sc nlib es sol}
\end{subfigure}
\hfill
\begin{subfigure}[b]{0.23\textwidth}
\centering
\includegraphics[width=\textwidth]{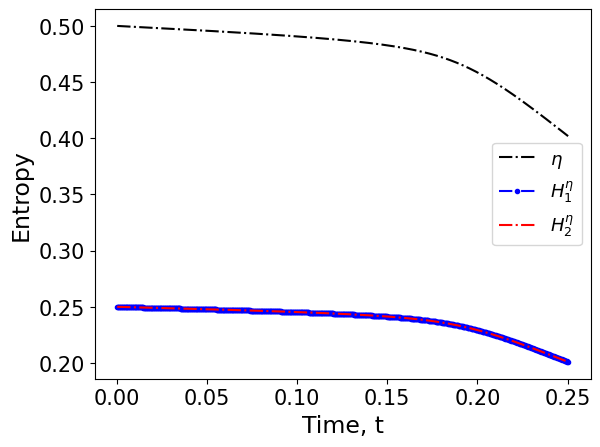}
\caption{Entropy functions}
\label{sc nlib es ent}
\end{subfigure}
\hfill 
\begin{subfigure}[b]{0.23\textwidth}
\centering
\includegraphics[width=\textwidth]{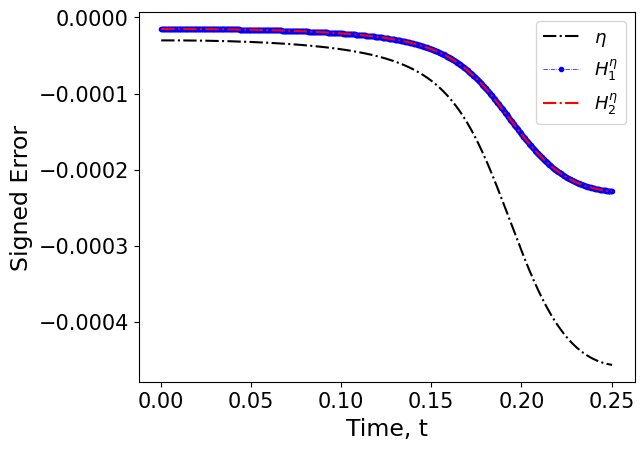}
\caption{Signed errors}
\label{sc nlib es Serr}
\end{subfigure}
\hfill
\begin{subfigure}[b]{0.23\textwidth}
\centering
\includegraphics[width=\textwidth]{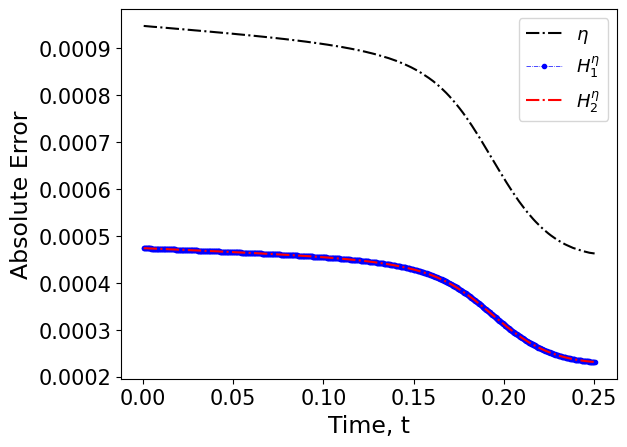}
\caption{Absolute errors}
\label{sc nlib es USerr}
\end{subfigure}
\caption{\centering Non-linear inviscid Burgers' test at $T=0.25$ using first order ES scheme with $C=0.1$ and $Nx=256$}
\label{sc nlib es}
\end{figure}

For this non-linear one-dimensional problem with $G^{(1)}(U)=\frac{1}{2}U^2$, we choose $\eta(U)=U^2$, and correspondingly $\omega^{(1)}(U)=\frac{2}{3}U^3$ satisfies the compatibility condition in \cref{Ent codn for mac model}. The initial condition is $U_0\left(x_1\right)=sin (2\pi x_1)$. The domain of the problem is $[0,1)$, and it is discretised using 256 uniform cells. Periodic boundary conditions are used here. We use entropy conserving and entropy stable schemes respectively for obtaining numerical solutions at $T=\frac{0.1}{2\pi} \text{ and } T=0.25$ in \cref{sc nlib,sc nlib es}.\\ 
\Cref{sc nlib sol,sc nlib es sol} show that the numerical solutions match well with the exact solutions. \Cref{sc nlib ent,sc nlib es ent} show that macroscopic and vector-kinetic entropy functions are conserved and dissipated respectively in the smooth ($T=\frac{0.1}{2\pi}$) and non-smooth ($T=0.25$) cases. The signed and absolute errors for $T=\frac{0.1}{2\pi}$ are shown in \cref{sc nlib Serr,sc nlib USerr}. Since we use second order accurate entropy conserving scheme for vector-kinetic model and $\Delta x$ is of $O(10^{-3})$, we expect an absolute error of $O(10^{-6})$ in the vector-kinetic entropies. However, we observe an absolute error of $O(10^{-4})$ in \cref{sc la USerr}. This might be because the terms multiplying $O(\Delta x^2)$ in the M-PDE of entropy equality are not $O(1)$ due to non-linearities. The negative signed errors in \cref{sc la Serr} indicate that the error is globally dissipative in nature. Due to the symmetric nature of periodic profile, there may be cancellations in errors spatially and we observe a very low signed error of $O(10^{-13})$. \\
Further, the signed and absolute errors for $T=0.25$ are shown in \cref{sc nlib es Serr,sc nlib es USerr}. Here too, we observe an absolute error of $O(10^{-4})$. Negative signed error of $O(10^{-4})$ indicates entropy dissipation after the formation of discontinuity. \\
In order to study the convergence of the problem, a very low CFL of $C=0.1$ is chosen. The reference solution is the exact solution obtained by employing Newton-Raphson iteration with tolerance of $10^{-15}$. It is seen from \cref{tab:EOC_la} that more than second order accuracy is attained as the grid is refined.  

\begin{table}[tbhp]
\begin{center}
\begin{tabular}{|m{2cm}|m{2cm}|m{2cm}|m{2cm}|}
\hline
\centering Number of cells, Nx & \centering $\Delta x_1$ & \centering $L_2$ norm & $O(L_2)$ \\
\hline
\centering 64 & \centering 0.015625 & \centering 0.000281831 &  - \\
\centering 128 & \centering 0.0078125 & \centering 0.000118395 &  1.89 \\
\centering 256 & \centering 0.00390625 & \centering 4.37E-05 &  3.24 \\
\hline
\end{tabular}
\caption{\centering EOC for non-linear inviscid Burgers' test at $T=\frac{0.1}{2\pi}$ using EC scheme with $C=0.1$} 
\label{tab:EOC_nlib}
\end{center}
\end{table}

\subsection{Shallow water equations}
\label{SubSec: SW}
We consider the shallow water equations,  
\begin{equation}
\partial_t \begin{bmatrix} \rho \\ \rho u_j \end{bmatrix} + \partial_{x_d} \begin{bmatrix} \rho u_d \\ \rho u_j u_d + p \delta_{dj} \end{bmatrix} = \bold{0} \ ; \ p=\kappa \rho^2 \ ; \ j \in \{1,2,..,D \}
\end{equation}
with initial condition $\bold{U}(x_1,..,x_d,..,x_D,0)=\bold{U_0}(x_1,..,x_d,..,x_D)$. Here, $\bold{U}= \begin{bmatrix} \rho \\ \rho u_j \end{bmatrix}$, $\bold{G}^{(d)}\left( \bold{U} \right)=\begin{bmatrix} \rho u_d \\ \rho u_j u_d + p \delta_{dj} \end{bmatrix}$ and $\kappa=\frac{1}{2}$. The notation $h, \ g$ with $h=\rho$, $g=2\kappa=1$ is commonly used in the shallow water community. In this case, $p=\frac{1}{2}g h^2$. \\
The entropy function is $\eta\left( \bold{U}\right)=\frac{1}{2}\rho u_j u_j + \kappa \rho^2$, and correspondingly the entropy flux functions become $\omega^{(d)}\left( \bold{U}\right)=u_d \left( \frac{1}{2}\rho u_j u_j + 2 \kappa \rho^2 \right)$. $\bold{F}_m$ and $H^{\eta}_m$ of vector-kinetic model are found using \cref{F def vec kin model} and \cref{H def vec kin model} respectively. The constants $a_m, b^{(d)}_m$ and $\lambda$ are chosen as described in \ref{app constants a,b}. The time step is chosen as
\begin{equation}
\Delta t \leq C \frac{\Delta x }{\lambda} \ ; \  \Delta x = min\left( \Delta x_d \right)
\end{equation}
Here, $C$ is the CFL number. Let us construct the entropy conserving flux $\left(v^{(d)}_m\bold{F}_m\right)^{\star}_{i_{d}+\frac{1}{2}}$ satisfying \cref{EC condn for vec kin model}. Consider the arithmetic average $\overline{A}_{i_{d}+\frac{1}{2}}=\frac{1}{2} \left( A_i + A_{i_{d}+1} \right)$. This average satisfies $\left[\left[AB\right]\right]_{i_{d}+\frac{1}{2}}=\overline{A}_{i_{d}+\frac{1}{2}}\left[\left[B\right]\right]_{i_{d}+\frac{1}{2}}+\overline{B}_{i_{d}+\frac{1}{2}}\left[\left[A\right]\right]_{i_{d}+\frac{1}{2}}$. Hence, the entropy conserving condition in \cref{EC condn for vec kin model} can be expressed as,
\begin{multline}
\left< \begin{bmatrix} 2\kappa \left[\left[\rho \right]\right]_{i_{d}+\frac{1}{2}} - \overline{u_k}_{i_{d}+\frac{1}{2}}\left[\left[u_k\right]\right]_{i_{d}+\frac{1}{2}} \\ \left[\left[u_j\right]\right]_{i_{d}+\frac{1}{2}} \end{bmatrix}, \left(v^{(d)}_m\bold{F}_{m}\right)^{\star}_{i_{d}+\frac{1}{2}} \right> = \\ v^{(d)}_m \left( 2\overline{\rho}_{i_{d}+\frac{1}{2}}\left(a_m \left[\left[\rho\right]\right]_{i_{d}+\frac{1}{2}} + b^k_m \overline{u_k}_{i_{d}+\frac{1}{2}}\left[\left[\rho\right]\right]_{i_{d}+\frac{1}{2}}\right) + \overline{\rho^2}_{i_{d}+\frac{1}{2}} \left( b^k_m \left[\left[u_k\right]\right]_{i_{d}+\frac{1}{2}} \right) \right) 
\end{multline}
Equating the terms corresponding to $\left[\left[\rho\right]\right]_{i_{d}+\frac{1}{2}}$ and $\left[\left[u_j\right]\right]_{i_{d}+\frac{1}{2}}$, we obtain
\begin{equation}
\left(v^{(d)}_m\bold{F}_{m}\right)^{\star}_{i_{d}+\frac{1}{2}}= \begin{bmatrix} v^{(d)}_m \overline{\rho}_{i_{d}+\frac{1}{2}} \left(a_m + b^k_m \overline{u_k}_{i_{d}+\frac{1}{2}}\right) \\ v^{(d)}_m \left( \overline{\rho}_{i_{d}+\frac{1}{2}} \overline{u_j}_{i_{d}+\frac{1}{2}} \left(a_m + b^k_m \overline{u_k}_{i_{d}+\frac{1}{2}}\right) +\kappa b^j_m\overline{\rho^2}_{i_{d}+\frac{1}{2}}  \right) \end{bmatrix}
\end{equation}
This EC flux is second order accurate in space. Let us now derive the entropy stable flux given by \cref{ES int flux for vec kin model}. We know that $\sum_{m=1}^M \bold{D}^{(d)}_{m_{i_{d}+\frac{1}{2}}}=\bold{D}^{(d)}_{i_{d}+\frac{1}{2}}$, a positive-definite matrix. We use the robust $\bold{D}^{(d)}_{i_{d}+\frac{1}{2}}$ described in \cite{doi:10.1137/110836961}. That is,
\begin{equation}
\bold{D}^{(d)}_{i_{d}+\frac{1}{2}}=\bold{R}^{(d)}_{i_{d}+\frac{1}{2}} \bold{\Lambda}^{(d)}_{i_{d}+\frac{1}{2}} \bold{R}^{(d)^T}_{i_{d}+\frac{1}{2}}
\end{equation}
where $\bold{R}^{(d)}_{i_{d}+\frac{1}{2}}$ is a suitably scaled matrix whose columns are eigenvectors of $\partial_{\bold{U}} \bold{G}^{(d)}$, and $\bold{\Lambda}^{(d)}_{i_{d}+\frac{1}{2}}$ is the Roe-type diffusion matrix (arithmetic averages are used). The matrices $\bold{R}^{(d)}_{i_{d}+\frac{1}{2}}$ and $\bold{\Lambda}^{(d)}_{i_{d}+\frac{1}{2}}$ for shallow water equations can be found in \cite{fjordholm_mishra_tadmor_2009}. Then, we use $\bold{D}^{(d)}_{m_{i_{d}+\frac{1}{2}}}=\frac{1}{M} \bold{D}^{(d)}_{i_{d}+\frac{1}{2}}, \ \forall m$, and these are positive-definite. \\
This results in a first order accurate ES flux. Let us derive the second order accurate ES flux given by \cref{HO ES int flux for vec kin model}. As in \cite{doi:10.1137/110836961}, we express $\bold{D}^{(d)}_{i_{d}+\frac{1}{2}} \left<\left< \bold{V}\right>\right>_{i_{d}+\frac{1}{2}} = \bold{R}^{(d)}_{i_{d}+\frac{1}{2}} \bold{\Lambda}^{(d)}_{i_{d}+\frac{1}{2}} \left<\left< \widetilde{\bold{W}}\right>\right>_{i_{d}+\frac{1}{2}}$ where $\left<\left< \widetilde{\bold{W}}\right>\right>_{i_{d}+\frac{1}{2}}=\bold{B}^{(d)}_{i_{d}+\frac{1}{2}} \bold{R}^{(d)^T}_{i_{d}+\frac{1}{2}} \left[\left[ \bold{V}\right]\right]_{i_{d}+\frac{1}{2}}$. Here, $\bold{B}^{(d)}_{i_{d}+\frac{1}{2}}$ is a positive diagonal matrix. Now, consider the minmod limiter
\begin{equation}
\mu (A,B) = \left\{ \begin{matrix} \text{s } min(|A|,|B|) & \text{if s}=sign(A)=sign(B) \\ 0 & \text{otherwise} \end{matrix}\right. 
\end{equation}
Then, the reconstruction
\begin{multline}
\left<\left< \widetilde{\bold{W}}\right>\right>_{i_{d}+\frac{1}{2}} = \bold{R}^{(d)^T}_{i_{d}+\frac{1}{2}} \left[\left[ \bold{V}\right]\right]_{i_{d}+\frac{1}{2}} - \frac{1}{2} \left( \mu \left( \bold{R}^{(d)^T}_{i_{d}+\frac{1}{2}} \left[\left[ \bold{V}\right]\right]_{i_{d}+\frac{1}{2}}, \bold{R}^{(d)^T}_{i_{d}+\frac{1}{2}} \left[\left[ \bold{V}\right]\right]_{i_{d+\frac{3}{2}}} \right) \right. \\ + \left. \mu \left( \bold{R}^{(d)^T}_{i_{d}+\frac{1}{2}} \left[\left[ \bold{V}\right]\right]_{i_{d}-\frac{1}{2}}, \bold{R}^{(d)^T}_{i_{d}+\frac{1}{2}} \left[\left[ \bold{V}\right]\right]_{i_{d}+\frac{1}{2}} \right)  \right)
\end{multline}
results in a second order accurate ES flux. Since $\bold{B}^{(d)}_{i_{d}+\frac{1}{2}}$ is a positive diagonal matrix, the sign property
\begin{equation}
sign\left( \left<\left< \widetilde{\bold{W}}\right>\right>_{i_{d}+\frac{1}{2}} \right) = sign\left( \bold{R}^{(d)^T}_{i_{d}+\frac{1}{2}} \left[\left[ \bold{V}\right]\right]_{i_{d}+\frac{1}{2}}  \right)
\end{equation}
holds true, and the entropy stability is maintained. For vector-kinetic entropy stability, we use $\bold{D}^{(d)}_{m_{i_{d}+\frac{1}{2}}} \left<\left< \bold{V}\right>\right>_{i_{d}+\frac{1}{2}} =\frac{1}{M} \bold{D}^{(d)}_{i_{d}+\frac{1}{2}} \left<\left< \bold{V}\right>\right>_{i_{d}+\frac{1}{2}}, \ \forall m$.  \\
It may be noted that we have derived the EC fluxes for vector-kinetic model from the vector-kinetic framework. Unlike this, we obtained the ES fluxes for vector-kinetic model based on the diffusion matrices commonly used in literature for macroscopic model. This is because the only requirement for entropy stability is positive-definiteness of $\bold{D}^{(d)}_{m_{i_{d}+\frac{1}{2}}}$, and we achieve this simply by employing the robust $\bold{D}^{(d)}_{i_{d}+\frac{1}{2}}$ used for macroscopic model. 

\subsubsection{1D expansion problem}
\begin{figure}
\centering
\begin{subfigure}[b]{0.26\textwidth}
\centering
\includegraphics[width=\textwidth]{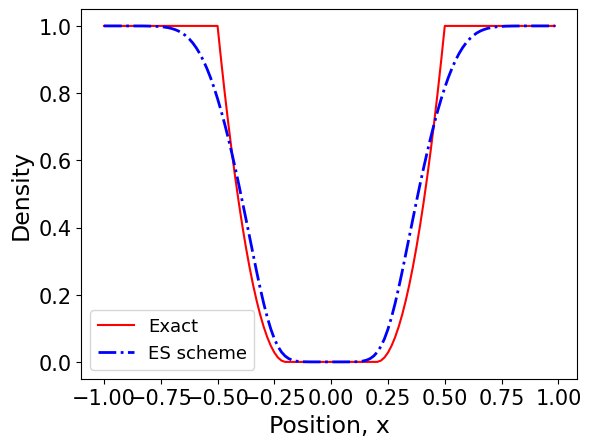}
\caption{$\rho$}
\label{sw exp sol_den}
\end{subfigure}
\hfill
\begin{subfigure}[b]{0.26\textwidth}
\centering
\includegraphics[width=\textwidth]{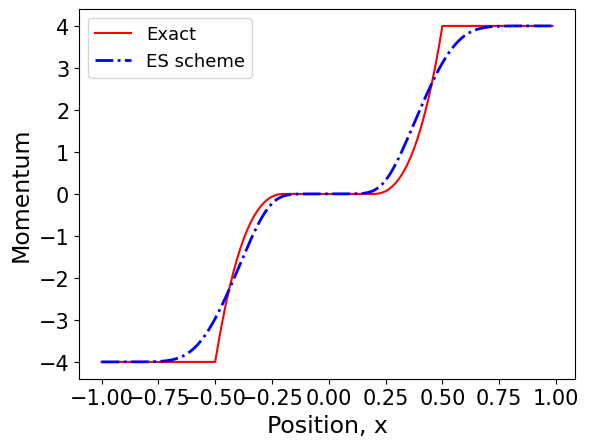}
\caption{$\rho u_1$}
\label{sw exp sol_momx}
\end{subfigure}
\hfill
\begin{subfigure}[b]{0.26\textwidth}
\centering
\includegraphics[width=\textwidth]{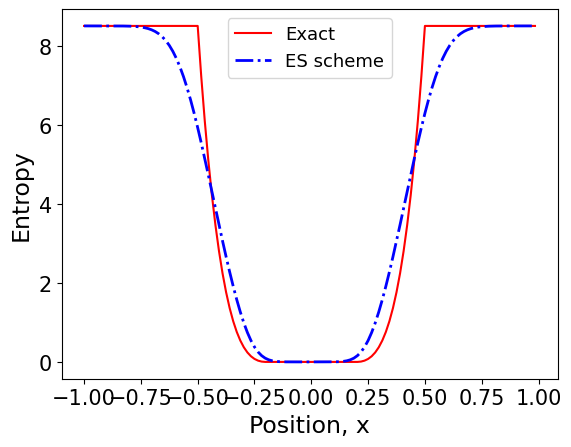}
\caption{$\eta$}
\label{sw exp sol_ent}
\end{subfigure}
\vfill
\begin{subfigure}[b]{0.23\textwidth}
\centering
\includegraphics[width=\textwidth]{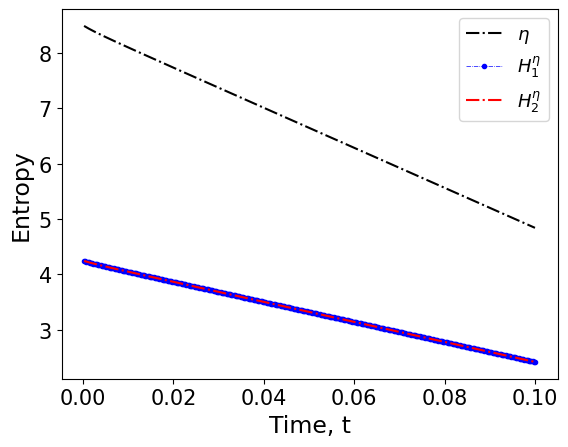}
\caption{Entropy functions}
\label{sw exp ent}
\end{subfigure}
\hfill
\begin{subfigure}[b]{0.26\textwidth}
\centering
\includegraphics[width=\textwidth]{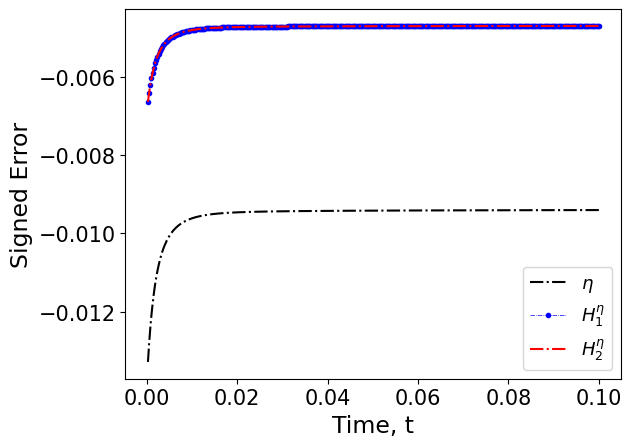}
\caption{Signed error}
\label{sw exp Serr}
\end{subfigure}
\hfill
\begin{subfigure}[b]{0.26\textwidth}
\centering
\includegraphics[width=\textwidth]{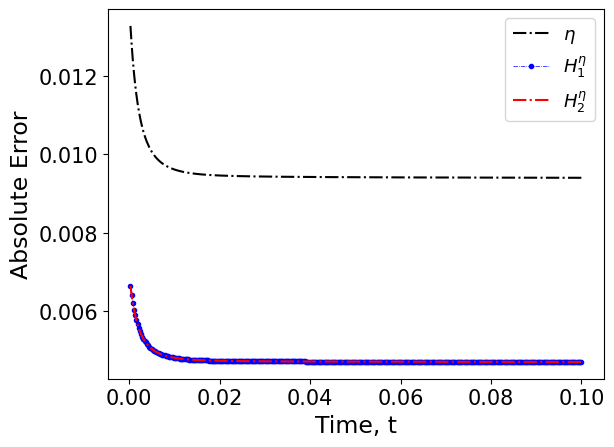}
\caption{Absolute error}
\label{sw exp USerr}
\end{subfigure}
\caption{\centering SW 1D expansion problem at $T=0.1$ using first order ES scheme with $C=0.1$ and $Nx=128$}
\label{sw exp}
\end{figure}

This test case is taken from \cite{fjordholm_mishra_tadmor_2009}. The domain of the problem is $[-1,1)$, and it is discretised using 128 uniform cells. The initial condition is,
\begin{equation}
\rho(x_1,0)=1, \ u_1(x_1,0)=\left\{  \begin{matrix} -4 & \text{if } x_1<0 \\ 4 & \text{if } x_1\geq0 \end{matrix} \right.
\end{equation}
Since the density can become very small, non-robust schemes will crash due to the in-ability to maintain positivity of density. Both entropy conserving and second order entropy stable schemes do not maintain the positivity. Hence, we utilise the first order entropy stable flux for vector-kinetic model to obtain the numerical results at $T=0.1$. The boundary values are kept fixed throughout the computation, and a very low CFL of $C=0.1$ is used for robustness. \\
It can be seen from \cref{sw exp sol_den} that the density remains non-negative. Further, the numerical solutions of density, momentum and entropy match well with the exact solution as shown in \cref{sw exp sol_den,sw exp sol_momx,sw exp sol_ent}. \Cref{sw exp ent,sw exp Serr,sw exp USerr} show entropy functions, their signed and absolute errors over time (for both macroscopic and vector-kinetic entropies). Since $\Delta x$ is of $O(10^{-2})$, one would expect an absolute error of $O(10^{-2})$ due to the usage of first order entropy stable flux. In \cref{sw exp USerr}, we observe a better absolute error of $O(10^{-3})$ in vector-kinetic entropies. Macroscopic entropy which is the sum of vector-kinetic entropies has an absolute error of $O(10^{-2})$. The negative signed errors in \cref{sw exp Serr} indicate the global dissipation of macroscopic and vector-kinetic entropies. This can also be seen in \cref{sw exp ent} from the decrease in global macroscopic and vector-kinetic entropies over time. It may be noted that the magnitudes of signed and absolute errors of all entropies in \cref{sw exp Serr,sw exp USerr} are same. This indicates that the first order entropy stable fluxes are dissipating the entropies at almost all spatial points, and not just globally. 

\subsubsection{1D dam break problem}
\begin{figure}
\centering
\begin{subfigure}[b]{0.26\textwidth}
\centering
\includegraphics[width=\textwidth]{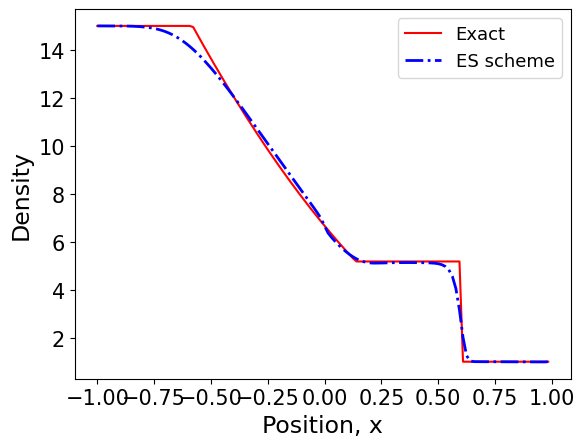}
\caption{$\rho$}
\label{sw db_1es sol_den}
\end{subfigure}
\hfill
\begin{subfigure}[b]{0.26\textwidth}
\centering
\includegraphics[width=\textwidth]{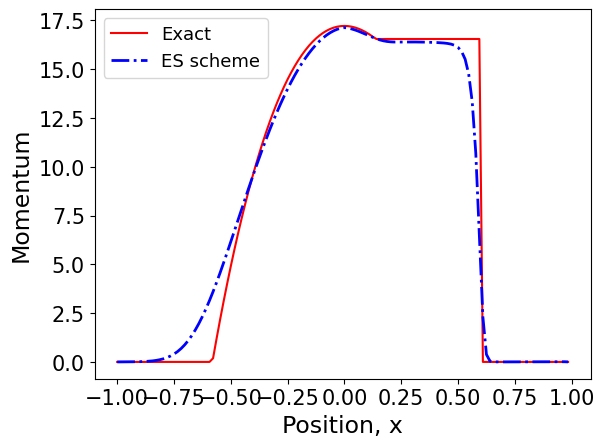}
\caption{$\rho u_1$}
\label{sw db_1es sol_momx}
\end{subfigure}
\hfill
\begin{subfigure}[b]{0.26\textwidth}
\centering
\includegraphics[width=\textwidth]{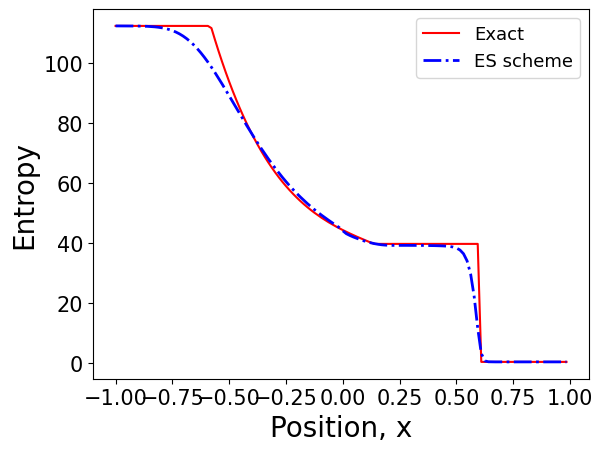}
\caption{$\eta$}
\label{sw db_1es sol_ent}
\end{subfigure}
\vfill
\begin{subfigure}[b]{0.26\textwidth}
\centering
\includegraphics[width=\textwidth]{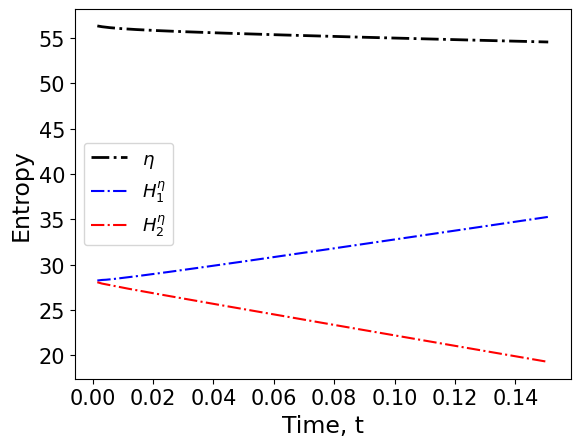}
\caption{Entropy functions}
\label{sw db_1es ent}
\end{subfigure}
\hfill 
\begin{subfigure}[b]{0.26\textwidth}
\centering
\includegraphics[width=\textwidth]{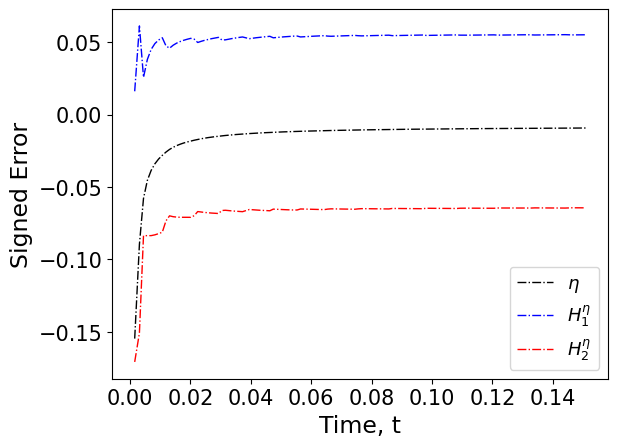}
\caption{Signed errors}
\label{sw db_1es Serr}
\end{subfigure}
\hfill
\begin{subfigure}[b]{0.26\textwidth}
\centering
\includegraphics[width=\textwidth]{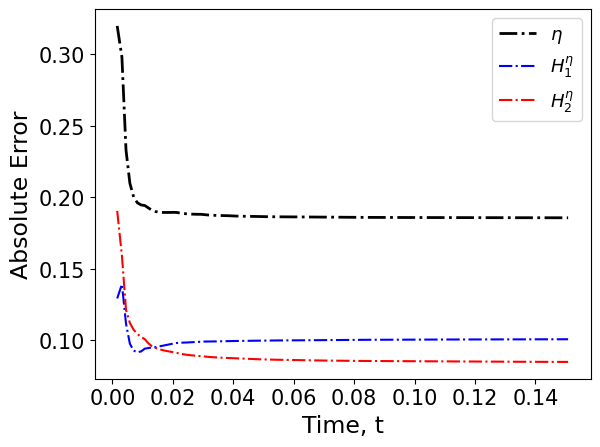}
\caption{Absolute errors}
\label{sw db_1es USerr}
\end{subfigure}
\caption{\centering SW 1D dambreak problem at $T=0.15$ using first order ES scheme with $C=0.4$ and $Nx=128$}
\label{sw db_1es}
\end{figure}

\begin{figure}[t]
\centering
\begin{subfigure}[b]{0.26\textwidth}
\centering
\includegraphics[width=\textwidth]{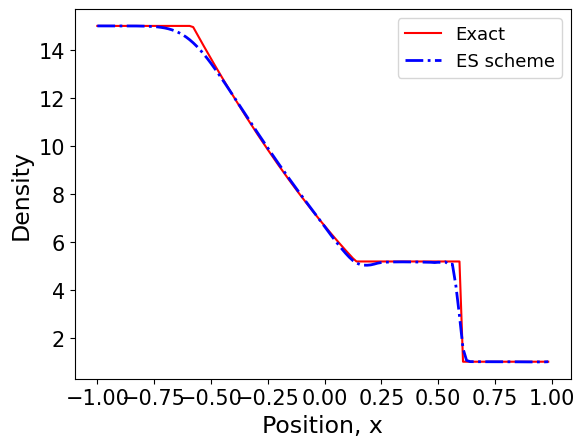}
\caption{$\rho$}
\label{sw db_2es sol_den}
\end{subfigure}
\hfill
\begin{subfigure}[b]{0.26\textwidth}
\centering
\includegraphics[width=\textwidth]{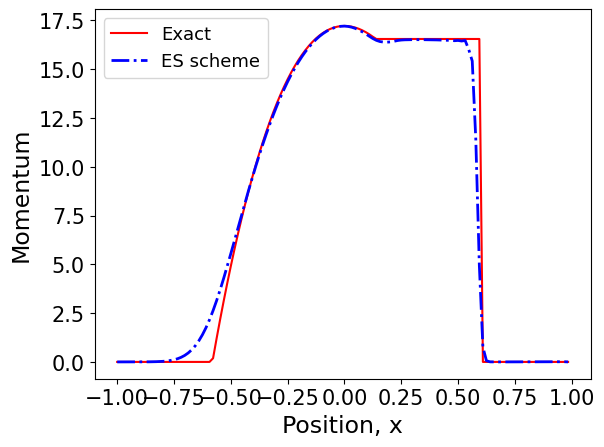}
\caption{$\rho u_1$}
\label{sw db_2es sol_momx}
\end{subfigure}
\hfill
\begin{subfigure}[b]{0.26\textwidth}
\centering
\includegraphics[width=\textwidth]{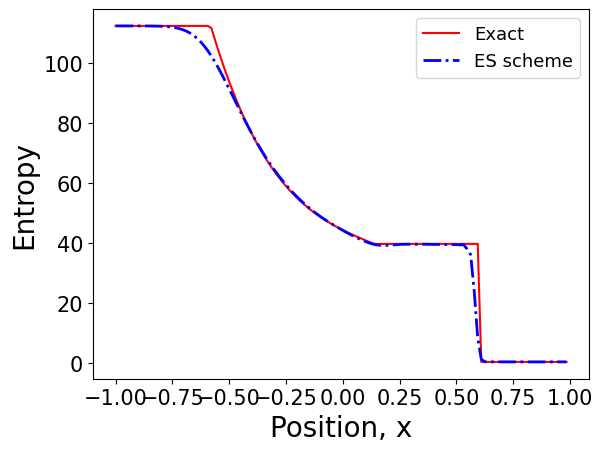}
\caption{$\eta$}
\label{sw db_2es sol_ent}
\end{subfigure}
\vfill
\begin{subfigure}[b]{0.26\textwidth}
\centering
\includegraphics[width=\textwidth]{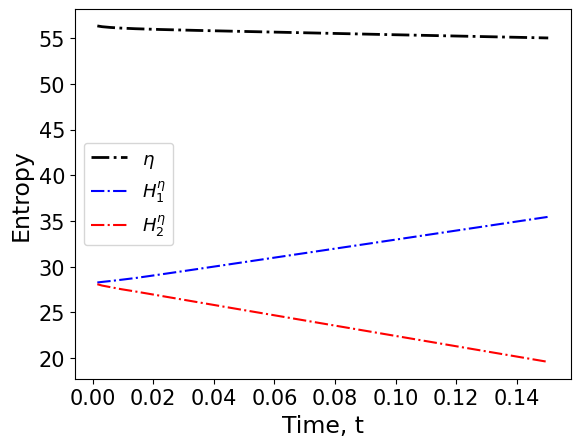}
\caption{Entropy functions}
\label{sw db_2es ent}
\end{subfigure}
\hfill 
\begin{subfigure}[b]{0.26\textwidth}
\centering
\includegraphics[width=\textwidth]{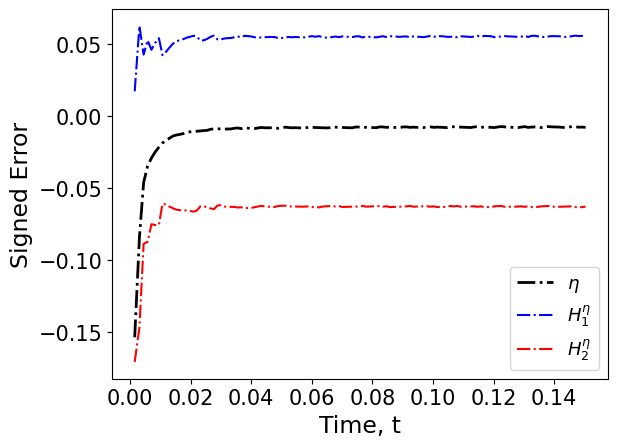}
\caption{Signed errors}
\label{sw db_2es Serr}
\end{subfigure}
\hfill
\begin{subfigure}[b]{0.26\textwidth}
\centering
\includegraphics[width=\textwidth]{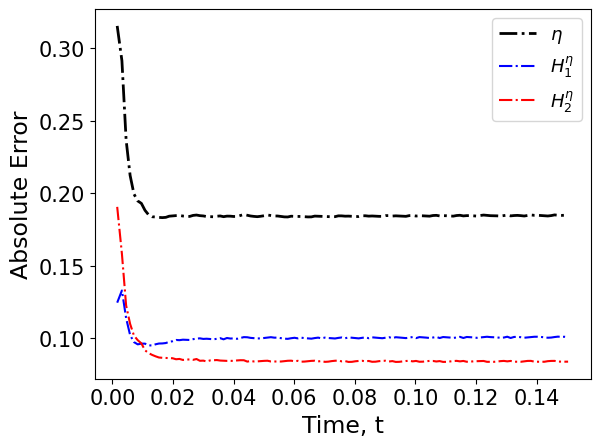}
\caption{Absolute errors}
\label{sw db_2es USerr}
\end{subfigure}
\caption{\centering SW 1D dambreak problem at $T=0.15$ using second order ES scheme (using minmod limiter) with $C=0.4$ and $Nx=128$}
\label{sw db_2es}
\end{figure}

This test case is also from \cite{fjordholm_mishra_tadmor_2009}. The domain of the problem is $[-1,1)$, and it is discretised using 128 uniform cells. The initial condition is,
\begin{equation}
\rho(x_1,0)=\left\{  \begin{matrix} 15 & \text{if } x_1<0 \\ 1 & \text{if } x_1\geq0 \end{matrix}, \ u_1(x_1,0)=0. \right.
\end{equation}
The numerical results obtained using first and second order entropy stable schemes at $T=0.15$ are shown in \cref{sw db_1es,sw db_2es} respectively. The second order entropy stable reconstruction need not produce monotone solutions near discontinuities. Hence, a minmod flux limiter (that combines first and second order entropy stable fluxes) is employed to produce monotone solution near discontinuities. The boundary values are kept fixed throughout the computation, and a CFL of $C=0.4$ is used. \\
It can be seen that both first and second order (with minmod limiter) schemes capture the solution profile reasonably well. A positive signed error for $H^{\eta}_1$ in \cref{sw db_1es Serr,sw db_2es Serr} indicates that the numerical diffusion added for the flux corresponding to $H^{\eta}_1$ is not sufficient to account for the entropy dissipation across discontinuities. This is because we have added equal weights of robust $\bold{D}^{(d)}_{i_{d}+\frac{1}{2}}$ to each of the vector-kinetic entropies, irrespective of their entropy dissipation requirements. Nevertheless, the error in macroscopic entropy which is obtained as the sum of vector-kinetic entropies is still negative (indicating entropy dissipation).

\subsubsection{2D periodic flow}
\begin{figure}[t]
\centering
\begin{subfigure}[b]{0.3\textwidth}
\centering
\includegraphics[width=\textwidth]{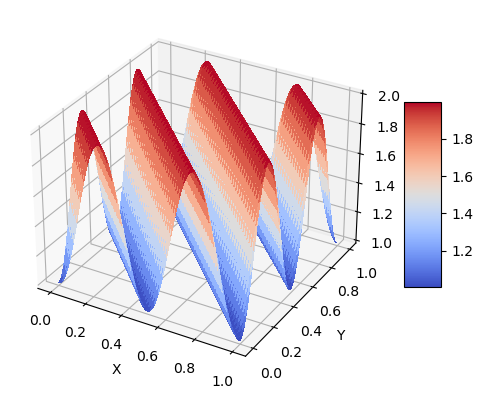}
\caption{Density at $T=0$}
\label{sw per ic}
\end{subfigure}
\hfill
\begin{subfigure}[b]{0.28\textwidth}
\centering
\includegraphics[width=\textwidth]{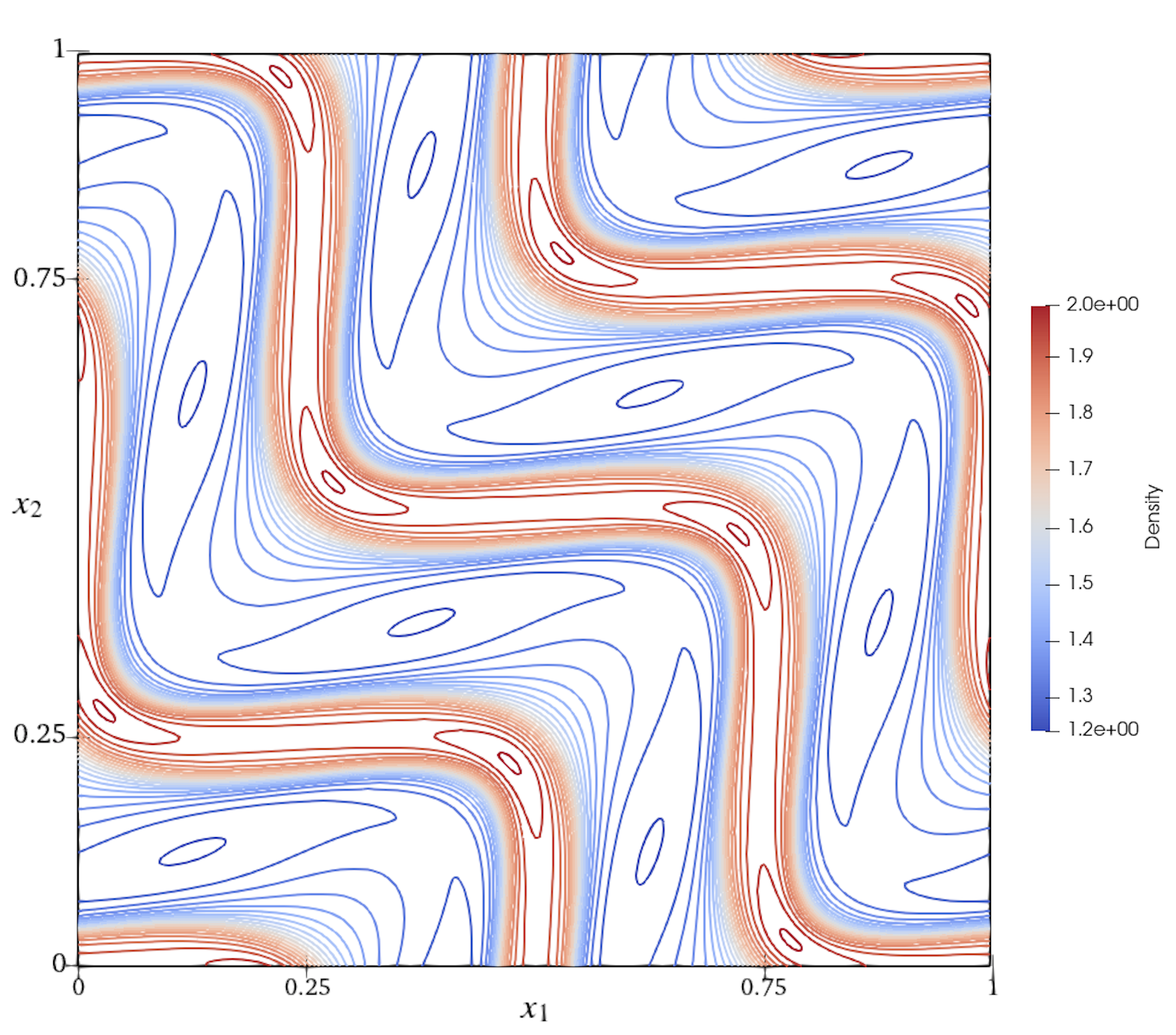}
\caption{Density contours at $T=0.1$}
\label{sw per sol}
\end{subfigure}
\vfill
\begin{subfigure}[b]{0.24\textwidth}
\centering
\includegraphics[width=\textwidth]{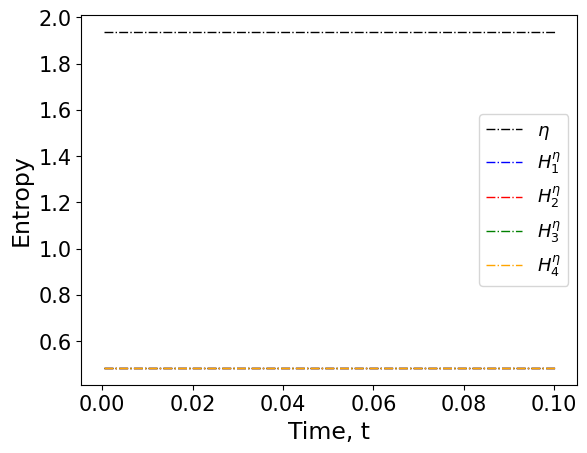}
\caption{Entropy functions}
\label{sw per ent}
\end{subfigure}
\hfill 
\begin{subfigure}[b]{0.235\textwidth}
\centering
\includegraphics[width=\textwidth]{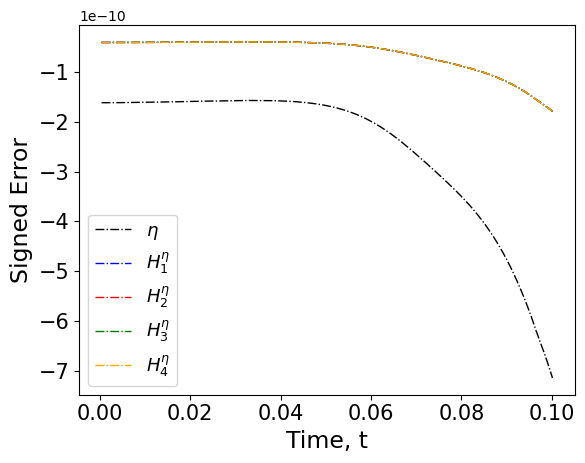}
\caption{Signed errors}
\label{sw per Serr}
\end{subfigure}
\hfill
\begin{subfigure}[b]{0.24\textwidth}
\centering
\includegraphics[width=\textwidth]{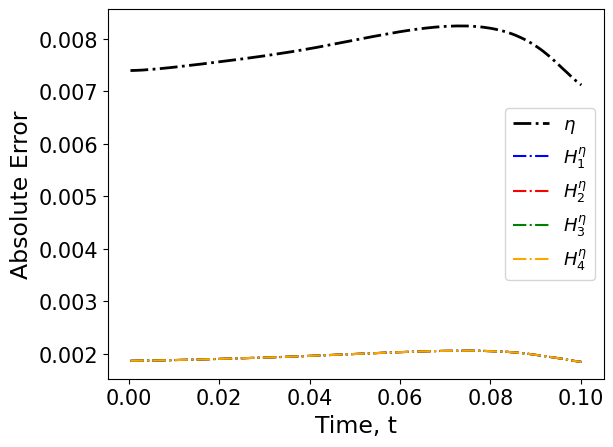}
\caption{Absolute errors}
\label{sw per USerr}
\end{subfigure}
\caption{\centering SW 2D periodic flow at $T=0.1$ using EC scheme with $C=0.5$ and $Nx, Ny=256$ (Blue, red and green lines are beneath the yellow line)}
\label{sw per}
\end{figure}

This test case is taken from the literature on asymptotic preserving schemes \cite{Kaiser2017}. In order to be useful in our context, we have taken the value of asymptotic parameter to be $1$. The domain of the problem is $[0,1)\times[0,1)$, and it is discretised using $256\times 256$ uniform cells. The initial condition shown in \cref{sw per ic} is given by,
\begin{eqnarray}
\rho(x_1,x_2,0)=1+sin^2 \left( 2 \pi \left( x_1+x_2\right)\right) \\ u_1(x_1,x_2,0)=u_2(x_1,x_2,0)=sin \left( 2\pi \left( x_1-x_2\right) \right)
\end{eqnarray}
The numerical results obtained using entropy conserving scheme at $T=0.1$ are shown in \cref{sw per sol}. Periodic boundary conditions are employed, and a CFL of C=0.5 is used. It can be seen from \cref{sw per ent} that the macroscopic and vector-kinetic entropy functions remain almost constant over time. From \cref{sw per USerr,sw per Serr}, we observe absolute and signed errors of $O(10^{-3})$ and $O(10^{-10})$ respectively. This huge difference implies that there are spatial cancellations between positive and negative errors. This may be due to the symmetric nature of periodic profile. Nevertheless, there is global dissipation of both macroscopic and vector-kinetic entropies as indicated by the negative errors in \cref{sw per Serr}. Order of convergence studies show that the accuracy attained is more than second order, and the results are shown in \cref{tab:EOC_sw per}. The reference solution for convergence studies is the numerical solution with refined grid of $512\times 512$.

\begin{table}[tbhp]
\begin{center}
\begin{tabular}{|m{0.5cm}|m{1.45cm}|m{1.4cm}|m{1cm}|m{1.4cm}|m{1.3cm}|m{1.4cm}|m{1.3cm}|}
\hline
\centering $N$ & \centering $\Delta x$ & \centering $||\rho||_{L_2}$ & $O(||\rho||)$ & \centering $||\rho u_1||_{L_2}$ & $O(||\rho u_1||)$ & \centering $||\rho u_2||_{L_2} $& $O(||\rho u_2||)$\\
\hline
\centering 32 & \centering 0.03125 & \centering 0.00162 &  - & \centering 0.00255 &  - & \centering 0.00255 &  - \\
\centering 64 & \centering 0.015625 & \centering 0.000378 &  2.10 & \centering 0.000362 &  2.82 & \centering 0.000362 &  2.82\\
\centering 128 & \centering 0.0078125 & \centering $5.64\times10^{-5}$ &  2.74 & \centering $5.54\times10^{-5}$ &  2.71  & \centering $5.54\times10^{-5}$ &  2.71 \\
\centering 256 & \centering 0.00390625 & \centering $7.62\times10^{-6}$ &  2.89 & \centering $7.33\times10^{-6}$ &  2.92 & \centering $7.33\times10^{-6}$ &  2.92 \\
\hline
\end{tabular}
\caption{\centering EOC for 2D periodic flow at $T=0.1$ using EC scheme with $C=0.5$} 
\label{tab:EOC_sw per}
\end{center}
\end{table}

\subsubsection{2D Travelling vortex}
\begin{figure}[t] 
\centering
\begin{subfigure}[b]{0.23\textwidth}
\centering
\includegraphics[width=\textwidth]{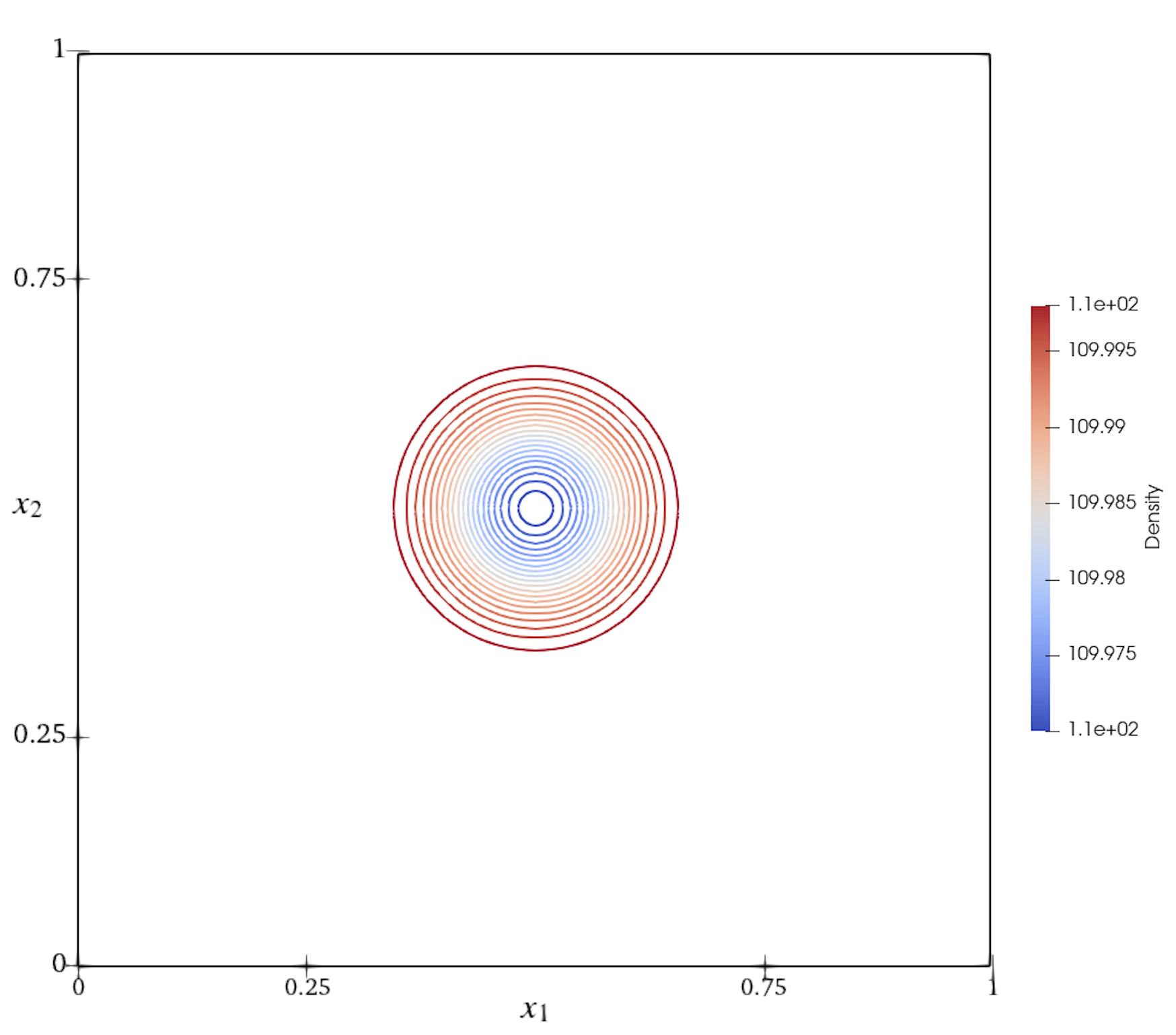}
\caption{Initial condn. $T=0$}
\label{sw tv_ec ic}
\end{subfigure}
\hfill
\begin{subfigure}[b]{0.23\textwidth}
\centering
\includegraphics[width=\textwidth]{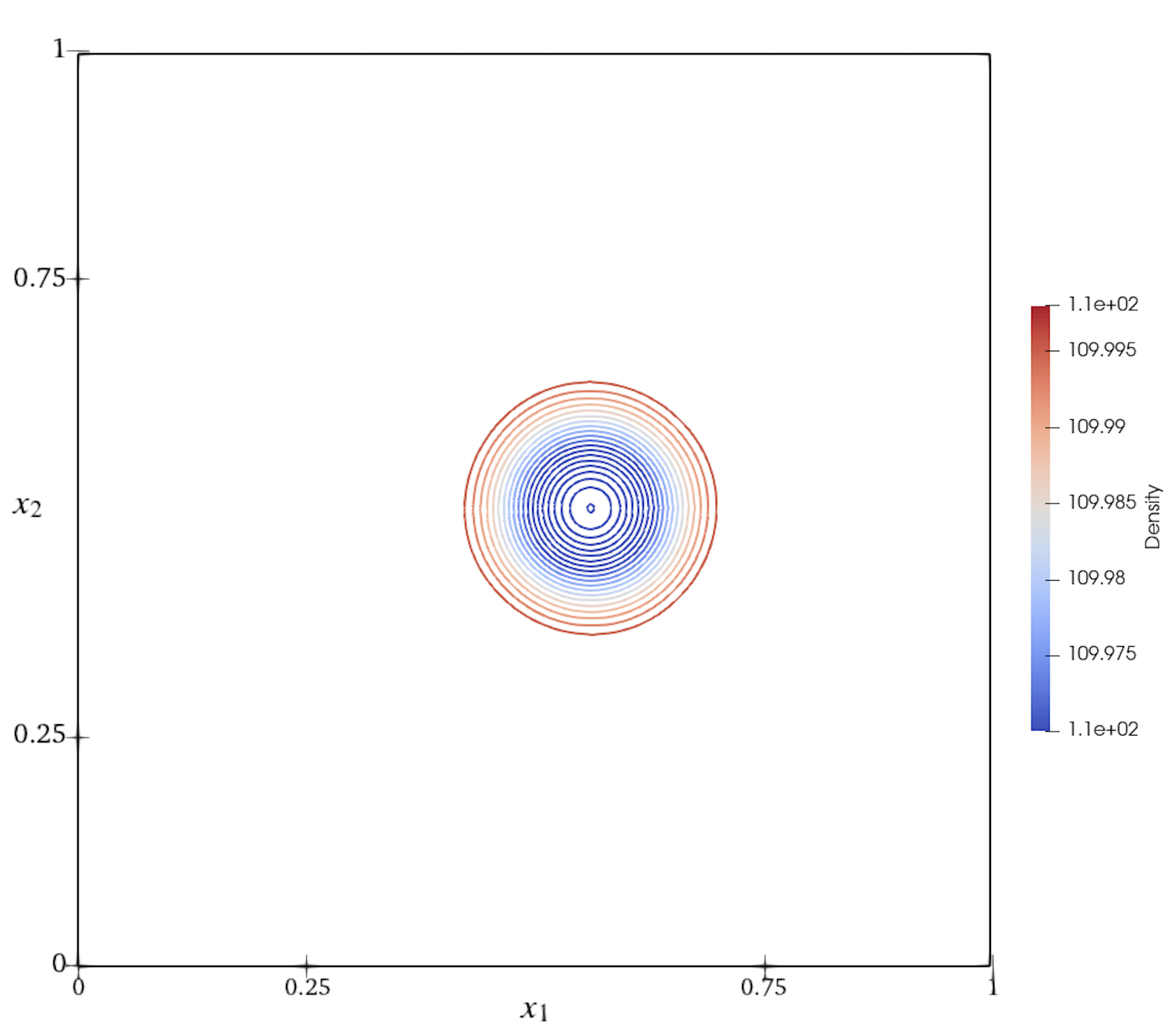}
\caption{Numerical soln. $T=0.1$}
\label{sw tv_ec sol}
\end{subfigure}
\hfill
\begin{subfigure}[b]{0.23\textwidth}
\centering
\includegraphics[width=\textwidth]{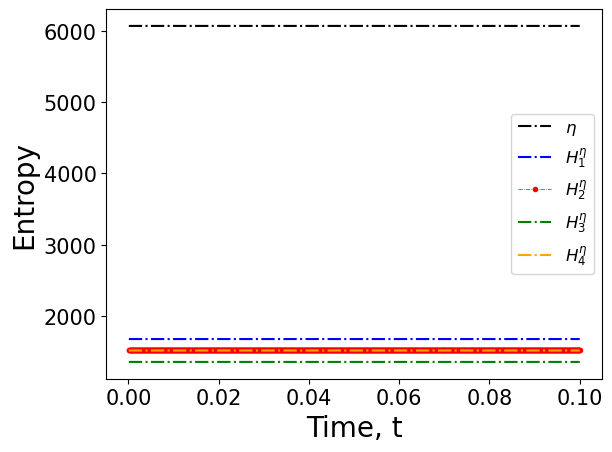}
\caption{Entropy functions}
\label{sw tv_ec ent}
\end{subfigure}
\hfill
\begin{subfigure}[b]{0.24\textwidth}
\centering
\includegraphics[width=\textwidth]{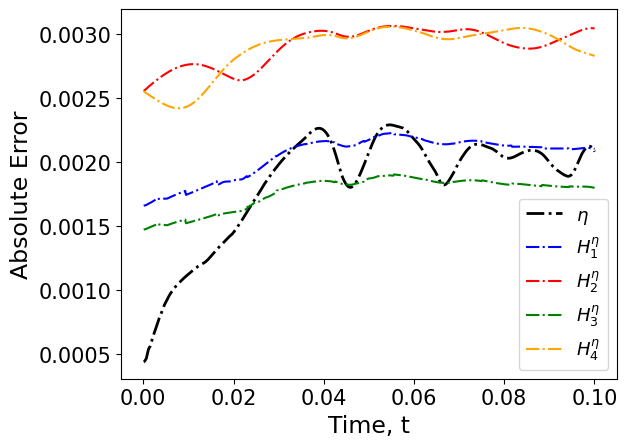}
\caption{Absolute errors}
\label{sw tv_ec USerr}
\end{subfigure}
\vfill
\begin{subfigure}[b]{0.26\textwidth}
\centering
\includegraphics[width=\textwidth]{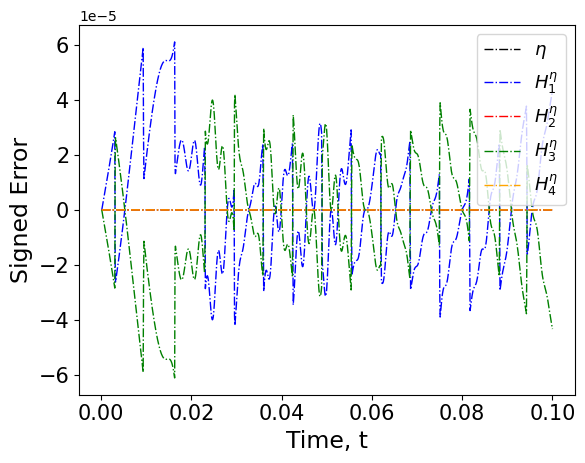}
\caption{Signed errors}
\label{sw tv_ec Serr}
\end{subfigure}
\hfill
\begin{subfigure}[b]{0.26\textwidth}
\centering
\includegraphics[width=\textwidth]{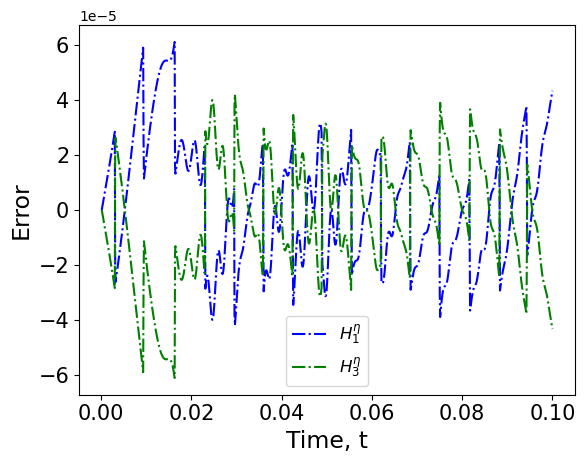}
\caption{$H^{\eta}_1$ and $H^{\eta}_3$}
\label{sw tv_ec H1H3 Serr}
\end{subfigure}
\hfill
\begin{subfigure}[b]{0.26\textwidth}
\centering
\includegraphics[width=\textwidth]{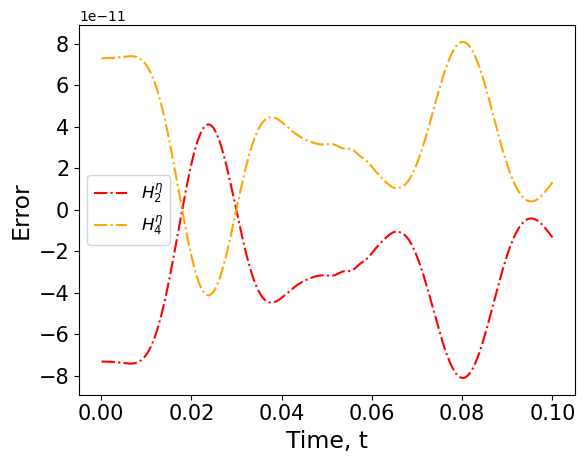}
\caption{$H^{\eta}_2$ and $H^{\eta}_4$}
\label{sw tv_ec H2H4 Serr}
\end{subfigure}
\caption{\centering SW 2D travelling vortex at $T=0.1$ using EC scheme with $C=0.5$ and $Nx, Ny=256$}
\label{sw tv_ec}
\end{figure}

This test case is also taken from the literature on asymptotic preserving schemes \cite{Kaiser2017}. We have taken the value of asymptotic parameter to be $0.8$, so that it will be useful in our context. The domain of the problem is $[0,1)\times[0,1)$, and it is discretised using $256\times 256$ uniform cells. The initial condition shown in \cref{sw tv_ec ic} is given by,
\begin{eqnarray}
\rho \left( x_1,x_2, 0\right)= 110 + \left( 0.64\left( \frac{1.5}{4\pi}\right)^2 \right) Drc\left(x_1,x_2\right) \left( k\left(rc\right)-k\left( \pi\right)\right) \\
u_1 \left( x_1,x_2, 0\right)=  0.6+ 1.5 \left( 1+cos \left( rc \left(x_1,x_2\right)\right)\right) Drc\left(x_1,x_2\right) \left( 0.5-x_2\right)  \\
u_2 \left( x_1,x_2, 0\right)=  0+ 1.5 \left( 1+cos \left( rc \left(x_1,x_2\right)\right)\right) Drc\left(x_1,x_2\right) \left( x_1-0.5\right) 
\end{eqnarray}
with 
\begin{eqnarray}
k\left(q\right)=2cos\left(q\right)+2q \ sin\left(q\right) + \frac{1}{8}cos\left(2q\right) + \frac{1}{4}q \ sin\left(2q\right) + \frac{3}{4}q^2 \\
rc \left( x_1,x_2\right) = 4 \pi \left(  \left( x_1-0.5\right)^2 + \left( x_2-0.5\right)^2 \right)^{\frac{1}{2}} \\
Drc\left(x_1,x_2\right) = \left\{ \begin{matrix}  1 & \text{if } rc\left( x_1,x_2\right)<\pi \\ 0 & \text{otherwise}\end{matrix} \right.
\end{eqnarray}
The second order entropy conserving and entropy stable schemes do not distort the structure of vortex, while the first order entropy stable scheme does. We present the numerical results obtained using second order entropy conserving scheme at $T=0.1$ as shown in \cref{sw tv_ec sol}. Periodic boundary conditions are employed, and a CFL of C=0.5 is used. \\
From \cref{sw tv_ec USerr}, we observe that the absolute errors of macroscopic and vector-kinetic entropies are of $O(10^{-3})$. On the other hand, the signed errors in $H^{\eta}_2$ and $H^{\eta}_4$ are of $O(10^{-11})$ (\cref{sw tv_ec H2H4 Serr}), while those in $H^{\eta}_1$ and $H^{\eta}_3$ are of $O(10^{-5})$ (\cref{sw tv_ec H1H3 Serr}). Moreover, the signed error profiles of vector-kinetic entropies are symmetric resulting in a much lower signed error of $O(10^{-14})$ for $\eta$ (not shown in plot). However, these symmetries in signed errors must be located at different spatial points. If they were located at the same spatial points, then we would observe a much lower absolute error in macroscopic entropy, unlike $O(10^{-3})$ in \cref{sw tv_ec USerr}. \\
Order of convergence studies are shown in \cref{tab:EOC_sw tv}. It is seen that the accuracy attained is more than second order for $\rho u_1$ and $\rho u_2$. For $\rho$, the required order of accuracy is observed in coarser mesh rather than in fine mesh, and this matches the conclusion made in \cite{Ricchiuto2021} where the analyses concerning types of vortices (based on their regularity) and their usage for validation of orders of accuracy of numerical methods are discussed. 

\begin{table}[tbhp]
\begin{center}
\begin{tabular}{|m{0.5cm}|m{1.45cm}|m{1.4cm}|m{1cm}|m{1.4cm}|m{1.3cm}|m{1.4cm}|m{1.3cm}|}
\hline
\centering $N$ & \centering $\Delta x$ & \centering $||\rho||_{L_2}$ & $O(||\rho||)$ & \centering $||\rho u_1||_{L_2}$ & $O(||\rho u_1||)$ & \centering $||\rho u_2||_{L_2} $& $O(||\rho u_2||)$\\
\hline
\centering 32 & \centering 0.03125 & \centering 0.000156 &  - & \centering 0.00339 &  - & \centering 0.00709 &  - \\
\centering 64 & \centering 0.015625 & \centering $4.39\times10^{-5}$ &  1.83 & \centering 0.000505 &  2.75 & \centering 0.00105 &  2.75\\
\centering 128 & \centering 0.0078125 & \centering $2.033\times10^{-5}$ &  1.11 & \centering $0.000105$ &  2.26  & \centering $0.000174$ &  2.60 \\
\hline
\end{tabular}
\caption{\centering EOC for 2D travelling vortex at $T=0.1$ using EC scheme with $C=0.5$} 
\label{tab:EOC_sw tv}
\end{center}
\end{table}

\subsubsection{2D cylindrical dambreak}
\begin{figure}
\centering
\begin{subfigure}[b]{0.23\textwidth}
\centering
\includegraphics[width=\textwidth]{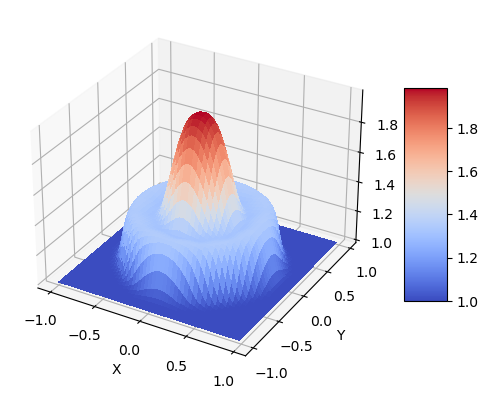}
\caption{Density at $T=0.2$}
\label{sw cdb_es1 sol}
\end{subfigure}
\hfill
\begin{subfigure}[b]{0.23\textwidth}
\centering
\includegraphics[width=\textwidth]{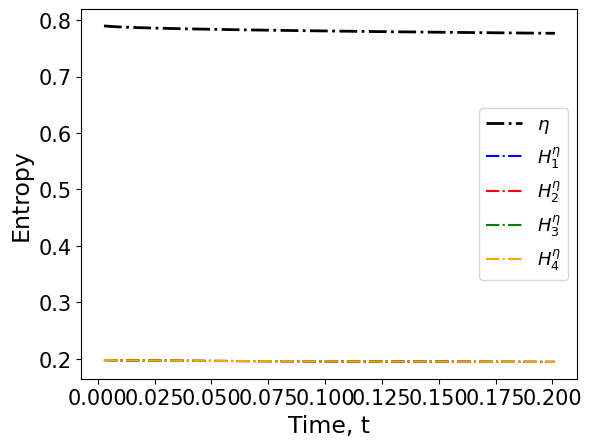}
\caption{Entropy functions}
\label{sw cdb_es1 ent}
\end{subfigure}
\hfill 
\begin{subfigure}[b]{0.23\textwidth}
\centering
\includegraphics[width=\textwidth]{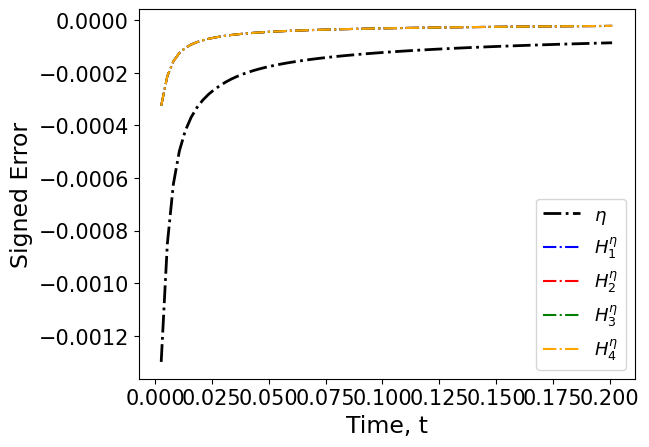}
\caption{Signed errors}
\label{sw cdb_es1 Serr}
\end{subfigure}
\hfill
\begin{subfigure}[b]{0.23\textwidth}
\centering
\includegraphics[width=\textwidth]{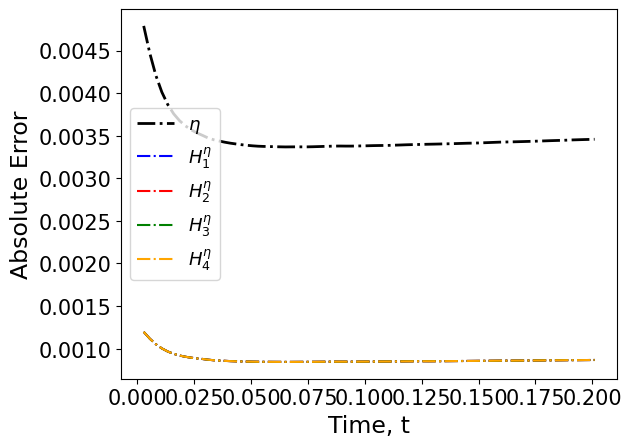}
\caption{Absolute errors}
\label{sw cdb_es1 USerr}
\end{subfigure}
\caption{\centering SW 2D cylindrical dam-break at $T=0.2$ using first order ES scheme with $C=0.4$ and $Nx, Ny=100$}
\label{sw cdb_es1}
\end{figure}

\begin{figure}
\centering
\begin{subfigure}[b]{0.23\textwidth}
\centering
\includegraphics[width=\textwidth]{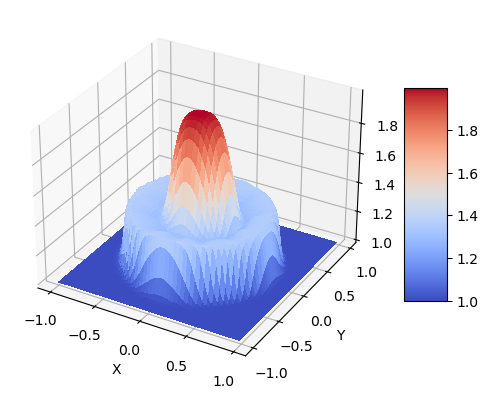}
\caption{Density at $T=0.2$}
\label{sw cdb_es2 sol}
\end{subfigure}
\hfill
\begin{subfigure}[b]{0.23\textwidth}
\centering
\includegraphics[width=\textwidth]{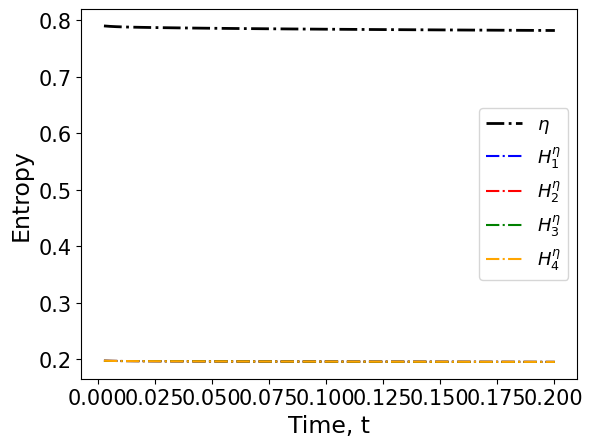}
\caption{Entropy functions}
\label{sw cdb_es2 ent}
\end{subfigure}
\hfill 
\begin{subfigure}[b]{0.23\textwidth}
\centering
\includegraphics[width=\textwidth]{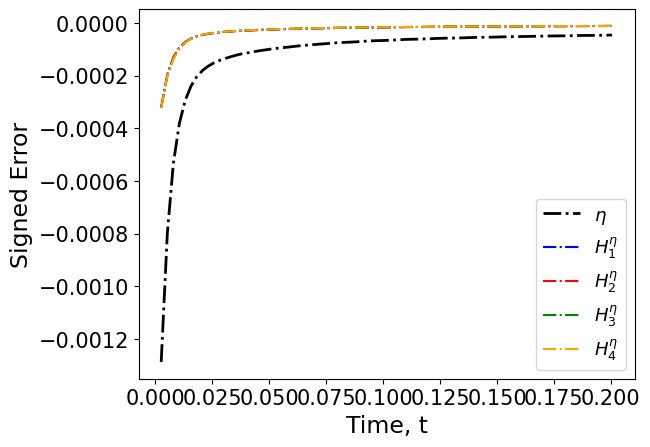}
\caption{Signed errors}
\label{sw cdb_es2 Serr}
\end{subfigure}
\hfill
\begin{subfigure}[b]{0.23\textwidth}
\centering
\includegraphics[width=\textwidth]{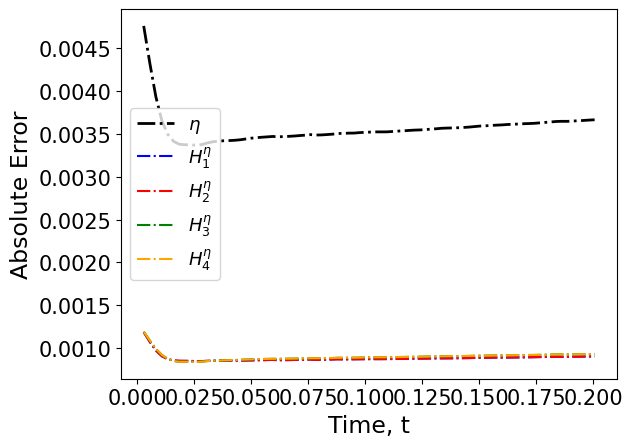}
\caption{Absolute errors}
\label{sw cdb_es2 USerr}
\end{subfigure}
\caption{\centering SW 2D cylindrical dambreak at $T=0.2$ using second order ES scheme (using minmod limiter) with $C=0.4$ and $Nx, Ny=100$}
\label{sw cdb_es2}
\end{figure}

This test case is taken from \cite{fjordholm_mishra_tadmor_2009}. The domain of the problem is $[-1,1)\times[-1,1)$, and it is discretised using $100\times 100$ uniform cells. The initial condition is given by,
\begin{equation}
\rho\left(x_1,x_2,0\right)= \left\{ \begin{matrix} 2 & \text{if } \left( x_1^2 + x_2^2 \right)^{\frac{1}{2}} < 0.5 \\ 1 & \text{otherwise}\end{matrix} \right. \ , \ u_1\left(x_1,x_2,0\right)=u_2\left(x_1,x_2,0\right)=0
\end{equation}
The numerical results of first and second order (with minmod limiter) entropy stable schemes at $T=0.2$ are shown in \cref{sw cdb_es1 sol,sw cdb_es2 sol} respectively. A CFL of $C=0.4$ is used, and periodic boundary conditions are employed. From \cref{sw cdb_es1 USerr,sw cdb_es2 USerr}, we observe that the absolute errors in entropies are of $O(10^{-3})$. Further, from \cref{sw cdb_es1 Serr,sw cdb_es2 Serr}, we observe that the signed errors in entropies are of $O(10^{-4})$. The negative signed errors indicate that there is global dissipation of entropy.  

\section{Summary and Conclusions}
\label{Sec:Conc}
 The following are the major highlights of the paper.
 \begin{itemize}
 \item{We provided a modification to the vector-BGK model, and this allows us to obtain entropy flux potentials that are required in the consistent definition of interface numerical entropy fluxes. \Cref{Exis ent var for vec kin model,pos def lem for vec kin model} are essential in obtaining the entropy flux potentials.}
 \item{We showed in \cref{thm EC for vec kin-mac model,thm ES for vec kin-mac model} that the moment of entropy conserving/stable schemes for vector-kinetic model results in entropy conserving/stable schemes for macroscopic model. \Cref{Exis ent var for vec kin model} plays a crucial role by rendering the linearities in the involved moments.}
 \item{In the numerical tests of scalar smooth problems, we employed our entropy conserving scheme and observed that the macroscopic and all the vector-kinetic entropies involved are conserved (up to absolute error). We also used signed error to observe global entropy dissipation/production due to higher order terms for which conservation does not apply.}
 \item{For shallow water equations, we derived an entropy conserving flux for vector-kinetic model by considering arithmetic averages of primitive variables. We used this entropy conserving scheme on smooth problems such as periodic flow and travelling vortex. In both cases, we observed the conservation of macroscopic and vector-kinetic entropies.}
 \item{We considered the 1D expansion problem where non-positivity of density can easily occur in non-robust schemes. For this, we employed the first order entropy stable scheme for vector-kinetic model and observed that the macroscopic and all vector-kinetic entropies involved are dissipative in nature. We also do not encounter non-positivity.}
 \item{In the non-smooth category, we considered scalar non-linear inviscid Burgers' test, 1D and 2D cylindrical dam-break problems. The second order entropy stable scheme employed for scalar case dissipates macroscopic and all vector-kinetic entropies. For the shallow water case, we employed the first and second order entropy stable schemes for vector-kinetic model. In 1D dam-break problem, we observed that some of the vector-kinetic entropies are not really dissipative, as their dissipation matrices are not built based on the dissipation requirements near discontinuities. Further research is required on the choice of appropriate robust dissipation matrices for vector-kinetic model. }
 \end{itemize} 
 Thus, the entropy preserving scheme developed in this paper preserves both vector-kinetic and macroscopic entropy functions. It is interesting to observe that the entropic numerical solutions of macroscopic model do not experience a notable difference when two different routes (via vector-kinetic and macroscopic) are taken. \\
  If the proposed entropy conserving scheme for vector-kinetic model is applied to the Euler's system, the vector-kinetic entropy conserving condition in \cref{EC condn for vec kin model} can be satisfied analogous to the ways available in literature to satisfy entropy conserving condition for macroscopic model in \cref{EC condn for mac model}. One can derive the fluxes by utilising an elegant and non-costly route available in literature (for instance, by defining primitive variables, substituting for entropy variables and entropy flux potentials in terms of these primitive variables into \cref{EC condn for vec kin model}, and equating the coefficients of the jumps in the primitive variables, as introduced in \cite{ISMAIL20095410} for satisfaction of the condition in \cref{EC condn for mac model}), and this is a work in progress. It is expected that the moment of such entropy conserving flux functions for vector-kinetic model derived using a particular method (say, \cite{ISMAIL20095410}) will be an entropy conserving flux function for macroscopic model derived using the same method (\cite{ISMAIL20095410}). 

\begin{appendix}
\section{Choice of constants $a_m, b^{(d)}_m$}
\label{app constants a,b} 
We know that the moment of \cref{vec kin model} becomes the given hyperbolic system in \cref{Mac model}, if the constants $a_m, b^{(d)}_m$ in \cref{F def vec kin model} satisfy the moment constraints in \cref{Mom cons. vec kin model_1,Mom cons. vec kin model_2}. We also know that, if the convex entropy function for vector-kinetic model $\left(\cref{H def vec kin model}\right)$ is used, then the moment of \cref{Ent eq for vec kin model} becomes \cref{Ent ineq for mac model} with equality. Further, positivity of eigenvalues of $\partial_{\bold{U}} \bold{F}_m$ is an important requirement for obtaining the entropy flux potentials and the results of \cref{thm EC for vec kin-mac model,thm ES for vec kin-mac model}. Therefore, in order for the formulation to hold, the constants $a_m, b^{(d)}_m$ are required to satisfy \cref{Mom cons. vec kin model_1,Mom cons. vec kin model_2} along with the positivity of eigenvalues of $\partial_{\bold{U}} \bold{F}_m$. \\
For one dimensional hyperbolic systems, we consider two discrete velocities, \textit{i.e.,} $M=2$. Let 
\begin{eqnarray}
a_1=\frac{1}{2}, a_2=\frac{1}{2} \\
 b^{(1)}_1=\frac{1}{2\lambda}, b^{(1)}_2=-\frac{1}{2\lambda}
 \end{eqnarray}
If $v^{(1)}_1=\lambda$ and $v^{(1)}_2=-\lambda$, then the moment constraints in \cref{Mom cons. vec kin model_1,Mom cons. vec kin model_2} are satisfied. Further, 
\begin{eqnarray}
\text{eig}\left(\partial_{\bold{U}} \bold{F}_1\right) = \text{eig}\left( \frac{1}{2} \bold{I} + \frac{1}{2\lambda} \partial_{\bold{U}} \bold{G}^{(1)}\right) \\
\text{eig}\left(\partial_{\bold{U}} \bold{F}_2\right) = \text{eig}\left(\frac{1}{2} \bold{I} -\frac{1}{2\lambda} \partial_{\bold{U}} \bold{G}^{(1)}\right)
\end{eqnarray} 
Thus, eigenvalues of $\partial_{\bold{U}} \bold{F}_m$ are $\frac{1}{2}\pm\frac{1}{2\lambda}\text{eig}\left( \partial_{\bold{U}} \bold{G}^{(1)}\right)$. Therefore, for positivity, we require $\lambda > \text{sup}\left( \left| \text{eig}\left( \partial_{\bold{U}} \bold{G}^{(1)} \right) \right| \right)$. The supremum is taken over all grid points/cells in the computational domain. \\
For two dimensional systems, we consider four discrete velocities, \textit{i.e.,} $M=4$. Let 
\begin{eqnarray}
a_1=\frac{1}{4}, a_2=\frac{1}{4}, a_3=\frac{1}{4}, a_4=\frac{1}{4}\\
b^{(1)}_1=\frac{1}{2\lambda}, b^{(1)}_2=0, b^{(1)}_3=-\frac{1}{2\lambda}, b^{(1)}_4=0\\
b^{(2)}_1=0, b^{(2)}_2=\frac{1}{2\lambda}, b^{(2)}_3=0, b^{(2)}_4=-\frac{1}{2\lambda}
\end{eqnarray}
If the following holds,
\begin{eqnarray}
v^{(1)}_1=\lambda, v^{(1)}_2=0, v^{(1)}_3=-\lambda, v^{(1)}_4=0 \\
v^{(2)}_1=0, v^{(2)}_2=\lambda, v^{(2)}_3=0, v^{(2)}_4=-\lambda 
\end{eqnarray}
then the moment constraints in \cref{Mom cons. vec kin model_1,Mom cons. vec kin model_2} are satisfied. Further, 
\begin{eqnarray}
\text{eig}\left(\partial_{\bold{U}} \bold{F}_1\right) = \text{eig}\left( \frac{1}{4} \bold{I} + \frac{1}{2\lambda} \partial_{\bold{U}} \bold{G}^{(1)}\right) \\
\text{eig}\left(\partial_{\bold{U}} \bold{F}_2\right) = \text{eig}\left( \frac{1}{4} \bold{I} +\frac{1}{2\lambda} \partial_{\bold{U}} \bold{G}^{(2)}\right) \\
\text{eig}\left(\partial_{\bold{U}} \bold{F}_3\right) = \text{eig}\left( \frac{1}{4} \bold{I} - \frac{1}{2\lambda} \partial_{\bold{U}} \bold{G}^{(1)}\right) \\
\text{eig}\left(\partial_{\bold{U}} \bold{F}_4\right) = \text{eig}\left( \frac{1}{4} \bold{I} -\frac{1}{2\lambda} \partial_{\bold{U}} \bold{G}^{(2)}\right)
\end{eqnarray} 
Thus, eigenvalues of $\partial_{\bold{U}} \bold{F}_m$ are $\frac{1}{4}\pm\frac{1}{2\lambda}\text{eig}\left( \partial_{\bold{U}} \bold{G}^{(1)}\right)$ and $\frac{1}{4}\pm\frac{1}{2\lambda}\text{eig}\left( \partial_{\bold{U}} \bold{G}^{(2)}\right)$. Therefore, for positivity, we require $\lambda > 2 \text{ sup}\left( \left| \text{eig}\left( \partial_{\bold{U}} \bold{G}^{(1)} \right) \right|, \left| \text{eig}\left( \partial_{\bold{U}} \bold{G}^{(2)} \right) \right| \right)$. The supremum is taken over all grid points/cells in the domain.
\end{appendix}

\section*{CRediT author statement}
\noindent \textbf{Megala Anandan}: Conceptualization, Methodology, Formal analysis, Software, Validation, Investigation, Writing- Original draft, Reviewing and Editing. \\
\textbf{S. V. Raghurama Rao}: Conceptualization, Writing- Reviewing and Editing.

\bibliographystyle{acm}
\bibliography{references}  
\end{document}